\documentclass[16pt]{article}
\usepackage{graphicx}
\usepackage{amssymb}
\usepackage{color}

 \let\cal\mathcal

\usepackage{graphicx}
\usepackage[curve]{xypic}
\usepackage{tikz}

\newcommand\C{{\mathbb C}}
\newcommand\Z{{\mathbb Z}}
\newcommand\N{{\mathbb N}}

 \newtheorem{theorem}{Theorem}[subsection]
\newtheorem{proposition}[theorem]{Proposition}
\newtheorem{corollary}[theorem]{Corollary}
\newtheorem{lemma}[theorem]{Lemma}
\newtheorem{definition}[theorem]{Definition}
\newtheorem{example}[theorem]{Example}
\newtheorem{remark}[theorem]{Remark}

\usepackage{graphicx}
 
\begin{document}
\title{ The Topology of Surface Singularities }
\author{Francoise Michel}

\maketitle

\begin{abstract}

 We consider a  reduced complex surface germ $(X,p).$  We do not assume that  $X$ is normal at $p$, and so,    the singular locus  $(\Sigma ,p)$  of $(X,p)$ could  be one dimensional. This text  is devoted to the description of the topology of $(X,p).$ By the conic structure theorem (see \cite{mi68}),  $(X,p)$ is homeomorphic to the cone on its link  $L_X.$

  First of all, any  good resolution,  $ \rho  :  (Y, E_{Y}) \to  (X, 0),$  of   $(X,p)$  factorizes    by the normalization $\nu : (\bar X,\bar p) \to (X,0 )$ (see \cite{la} Thm. 3.14). This  is why we proceed in two steps. 
  
  1)   When  $(X,p)$ a normal germ of surface,  $p$ is an isolated singular point  and  the link  $L_X$ of $(X,p)$ is a  well defined three differentiable manifold. Using the  good minimal resolution  of  $(X,p)$, $L_X$ is given as the boundary of  a well defined  plumbing (see Section 2) which has  a negative definite intersection form (see \cite{h-n-k} and \cite{ne}).

  2) In Section 3 , we use  a suitably  general morphism,   $\pi :  (X,p) \to  (\C ^2, 0)$,  to describe  the topology of a surface germ  $(X,p)$  which has a 1-dimensional singular locus $(\Sigma ,p)$. We detail the quotient  morphism induced by the normalization   $\nu$
  on the link $L_{\bar X}$  of $ ( \bar X, \bar p)$ (see also Section 2 in  Luengo-Pichon  \cite{l-p}).

 In Section 4,  we detail the proof of  the existence of a good resolution, of a normal surface germ,  by  the Hirzebruch-Jung  method (Theorem 4.2.1).   With this method  a good resolution   is obtained via an  embedded resolution of  the  discriminant of $\pi$ (see \cite{hi}).  An example is given Section 6.

 An  appendix (Section 5), is devoted to the topological study   of lens spaces and to the description of     the minimal  resolution  of quasi-ordinary  singularities of surfaces. Section 5 achieves to make the proof of Theorem 4.2.1 self-contained.

\end{abstract}

{\bf \small Mathematics Subject Classifications (2000).}   {\small 14B05, 14J17, 32S15,32S45, 32S55, 57M45}.

\bigskip

{\bf Key words. } {\small Surface singularities, Resolution of singularities, Normalization, 3-dimensional Plumbed Manifold, Discriminant. }

 \bigskip
 \section{Introduction}

 Let $I$ be a reduced ideal in $\C \{z_1, \dots ,z_n \} $ such that the quotient algebra $A_X=\C \{z_1, \dots ,z_n \} / I$ is two-dimensional.
The zero locus,  at the origin  $0$ of  $\C ^{n}$, of a set of generators of $I$ is an analytic surface germ  embedded in $(\C^{n},0).$   Let $(X,0)$ be its intersection with  the  compact   ball $B_{\epsilon}^{2n}$ 
of radius a sufficiently small  $\epsilon$,  centered at the origin in $\C^{n}$,  and $L_X$ its  intersection with the boundary   $S_{\epsilon} ^{2n-1}$ of $B_{\epsilon}^{2n}.$  Let $\Sigma $ be the set of the singular points of $(X,0).$ 

As $I$ is reduced $\Sigma$ is empty when $(X,0)$ is smooth, it is equal to the origin when $0$ is an isolated singular point, it is a curve when the germ has a non-isolated singular locus (in particular we don't exclude the cases of reducible germs). 

 If $\Sigma $ is a curve,   $K_{\Sigma }=\Sigma \cap S_{\epsilon} ^{2n-1}$ is the  disjoint union of $r$ one-dimensional circle ($r$ being the number of irreducible components of $\Sigma$) embedded in $L_X$. We  say that $K_{\Sigma }$ is the link  of $\Sigma.$ By the conic structure theorem (see \cite{mi68}), for a sufficiently small $\epsilon$,  $(X,\Sigma, 0)$ is homeomorphic to the cone on  the pair $(L_X, K_{\Sigma }) $ and to the cone on $L_X$ when $\Sigma =\{ 0 \} $.

 On the other hand,  thanks to A.Durfee  in \cite{du} , the  homeomorphism class of $(X,\Sigma ,0)$ depends only on the isomorphism class of the algebra $A_X$ (i.e. is independent of  a sufficiently small $\epsilon$ and of the choice of the embedding in $(\C ^n,0)$).  
We  say that the analytic type of $(X,0)$ is given by the isomorphism class of $A_X$ and,   we say that its  topological type  is given by the homeomorphism  class of the pair $(X, 0)$ if  $0$ is an isolated singular point, and  by the homeomorphism  class of the triple  $(X,\Sigma, 0)$ if the singular locus $\Sigma$ is a curve. 

\begin{definition}  {\bf The  link of $(X,0)$}  is the homeomorphism class of   $L_X$ if  $0$ is an isolated singular point (in particular if ($X,0$) is normal at $0$), and  of the  homeomorphic class of the pair   $(L_X, K_{\Sigma })$ if the singular locus $\Sigma$ is a curve. 

 \end{definition}

This paper is devoted to  the  description of the link of $(X,0).$ 

\subsection{Good resolutions}

\begin{definition} 
 A morphism $ \rho  :  (Y, E_{Y}) \to  (X, 0)$ where $E_{Y}=\rho ^{-1}(0) $ is the exceptional divisor of $\rho, $ is a {\bf  good   resolution}   of $(X,0)$ if :
 \begin{enumerate}

 \item  $Y$ is a smooth complex  surface, 
 \item the total transform  of the singular locus $E^+_{Y}=\rho ^{-1}(\Sigma) $ is a normal crossing divisor with smooth irreducible components.
 \item the restriction of $\rho $ on $Y\setminus  E^+_{Y}$ is an isomorphism.
\end{enumerate}

 \end{definition}
 
 \begin{definition}
 
 Let $\rho: (Y,E_Y)\longrightarrow (X,0)$ be a good resolution of $(X,0)$. 
 
 The  {\bf dual graph associated to $\rho$}, denoted $G(Y)$,   is constructed as follows.  The vertices  of $G(Y)$ represent the irreducible components of $E_Y.$   When two irreducible components of  $E_Y$ intersect, we join their associated vertices by edges which number is equal to the number of intersection points.  A  dual graph is a {\bf bamboo } if the graph is homeomorphic to a segment and if all the vertices represent a rational curve.

  If $E_i$ is an irreducible component of $E_Y$, let us denote by  $e_i$  the self intersection  of $E_i$ in $Y$ and  by $g_i$ its genus.  To obtain the   {\bf weighted  dual graph associated to $\rho$}, denoted $G_w(Y)$, we weight $G(Y)$   as follows. A  vertex  associated to the    irreducible $E_i$ of $E_Y$  is weighted by  $(e_i)$ when $g_i=0$ and by  $(e_i,g_i)$ when  $g_i>0$.

  \end{definition}
  
  For example if $X=\{ (x,y,z)\in \C ^{3},  z^m= x^k y^l\}, $ where  $m,k$ and $l$  are  integers greater than two and prime to by to, Figure 1 describes  the shape of  the dual graph of the minimal  good resolution of $(X,0)$.

 \begin{figure}[h]
$$
\unitlength=0.40mm
\begin{picture}(40.00,20.00)
\thicklines

\put(-85,10){$(e_1)$}
\put(70,10){$(e_n)$}
\put(-100,0){\line(1,0){90}}
\put(10,0){\line(1,0){90}}
\put(-105,-2.3){$<$}
\put(95,-2.3){$>$}
\put(-80,0){\circle*{4}}
\put(-60,0){\circle*{4}}
\put(-40,0){\circle*{4}}
\put(-20,0){\circle*{4}}
\put(-5,0){\ldots\ldots }
\put(20,0){\circle*{4}}
\put(80,0){\circle*{4}}
\put(60,0){\circle*{4}}
\put(40,0){\circle*{4}}

\end{picture}
$$

\caption{  $G_w(Y)$ when $X=\{ (x,y,z)\in \C ^{3},  z^m= x^k y^l\}.$ Here $G(Y)$ is a bamboo. The arrows represent the strict transform of $\{xy=0\}$.  
In particular if $m=12, k=5$ and $l=11$  the graph has three vertices with $e_1=-3, e_2=-2, e_3=-3$ (see \cite{m-p-w-1} p. 759).
} 
\end{figure}
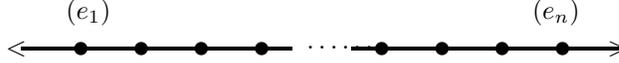

\begin{remark} 

If $(X,0)$ is reducible, let $(\cup _{1\leq i \leq r}X_i,0)$ be its decomposition as a union of irreducible surface  germs. Let $\nu_i : (\bar{X}_i,p_i) \to (X_i,0)$ be the normalization of the irreducible components of $(X,0)$. The  morphisms $\nu_i$ induce the   normalization morphism   on the disjoint union $\coprod _{1\leq i \leq r}(\bar{X}_i,p_i)$.

\end{remark}

\begin{remark}

  First of all, any  good resolution factorize   by the normalization $\nu : (\bar X,\bar p) \to (X,0 )$ (see \cite{la} Thm. 3.14).  In section 3, we  describe the topology of  normalization morphisms. After that it will be sufficient to describe the topology of the links of normal surface germs.

 A good resolution is minimal if its exceptional divisor doesn't contain any irreducible component of genus zero,  self-intersection $-1$ and which meets only one or two other irreducible components.  Let $ \rho  :  (Y, E_{Y}) \to  (X, 0)$ be a good resolution  and $\rho' : ({ Y'}, E_{ Y'})\to (X,0)$ be a good minimal resolution of $(X,0)$, then there exists   a morphism $\beta : (Y,E_Y)\to  (Y',E_{Y'})$  which is a  sequence of  blowing-downs of  irreducible components of genus zero and self-intersection $-1$ (see \cite{la} Thm 5.9 or \cite{b-p-v} p. 86). It implies the unicity, up to isomorphism, of the minimal good resolution of $(X,0).$

 As $\rho '$ factorize by $\nu$, the minimal good resolution of  $(\bar X,\bar p)$ can be  defined on $Y'.$ Let   $\bar {\rho }: ({  Y'}, E_{ Y'})\to (\bar X,\bar p)$  be the minimal good resolution of   $ (\bar X,\bar p)$ defined on $Y'.$  What we said just above implies that $\rho =\nu \circ \bar {\rho}\circ \beta $, i.e. $\rho$ is the composition of the three following morphisms:

$$(Y,E_Y) \stackrel {\beta }{\longrightarrow} (Y',E_{Y'}) \stackrel {\bar \rho }{\longrightarrow} (\bar X,\bar p) \stackrel {\nu }{\longrightarrow}  (X,0)$$

 \end{remark} 
 
 \subsection{Link of a complex surface germ}

  In Section 2, we describe the topology of a plumbing and  the topology of its  boundary. We explain how the existence of a good  resolution  describes  the link of a normal complex surface germ  as the boundary of a  plumbing  of disc bundles on oriented smooth compact real surfaces  with empty boundary. The boundary of a plumbing is, by definition,  a plumbed  3-manifold  (\cite{h-n-k}, \cite{ne}) or  equivalently a graph manifold in the sense of Waldhausen (\cite{wa}). The plumbing given by  the minimal  good resolution  of $(X,0)$ has a normal form in the sense of Neumann \cite {ne}  and represents its boundary in a unique way. 
 
 It implies that the link  of a  normal complex surface germ  $(X,0)$ determines the weighted dual  graph of its  good minimal resolution.   In particular,  if the link is $S^3$,  then  the good minimal resolution  of $(X,0)$ is an isomorphism and $(X,0)$ is smooth at the origin. This  is the famous result obtained in 1961 by Mumford \cite{mu}.  When the singular locus of $(X,0)$ is an irreducible germ of curve, its link can be $S^3$.     L\^e's conjecture, which is still open (see \cite{l-p} and \cite{debo} for partial results),  states that it can   only happen for  an equisingular family of irreducible curves.

 In  Section 3, we use a  suitably  general (as told   in 3.1)  projection  $\pi: (X,0) \to (\C^2,0)$  to describe  the topology of the restriction $\nu_L : L_{\bar X}  \to L_X$ of   the normalization  $\nu$ on the link $L_{\bar X}.$  We will show  that $\nu_L$ is a homeomorphism if and only if a general hyperplane section of $(X,0)$ is locally irreducible at $z$ for all points  $z\in (\Sigma \setminus \{ 0\}).$ Otherwise, as  stated without a  proof in  Luengo-Pichon  \cite{l-p},
 $\nu_L$ is  the composition of two kind of  topological quotients: curlings and identifications. Here, we give detailed proofs. Some years ago,   John Milnor asked me for a description of the link of a surface germ with non-isolated singular locus.  I hope that  Section 3 gives  a satisfactory answer.
 
 In Section 4 we suppose that $(X,0)$ is { \bf {normal}.} We use a finite morphism  $\pi: (X,0) \to (\C^2,0)$ and its discriminant $\Delta $,  to obtain a good resolution $ \rho  :  (Y, E_{Y}) \to  (X, 0)$  of $(X,0).$  We follow  Hirzebruch's method (see \cite{hi}, see also Brieskorn \cite{br2} for a presentation of Hirzebruch's work).  The scheme  to obtain $\rho$ is as in \cite{m-m}, but  our redaction here is quite different. In \cite{m-m},   the purpose is to  study   the behaviour of invariants  associated to finite morphisms defined   on $(X,0)$. Here, we detail the topology of each steps of the construction to precise the behaviour of  $\rho$.   Hirzebruch's  method uses the  properties of the topology of the normalization,  presented in Section 3,  and the resolution of the quasi-ordinary singularities of surfaces already studied by Jung. This  is why one says that this resolution $\rho$ is the Hirzebruch-Jung  resolution associated to $\pi.$   Then $L_X$ is homeomorphic to   the boundary of a regular neighborhood of the exceptional divisor $E_Y$ of  $ \rho  :  (Y, E_{Y}) \to  (X, 0)$ which is a plumbing as defined in Section 2. 
 
 Section 5  is  an appendix which can be read independently of the other sections.    We suppose  again that $(X,0)$ is { \bf {normal}.} We give topological proofs of basic results, already used in Section 4,    on finite morphism $\phi : (X,0) \to (\C^2,0)$,   in the two following cases:
 
 1) The discriminant of $\phi $ is a smooth germ of curve.  Then, in Lemma 5.2.1,   we show that $(X,0)$ is analytically isomorphic to $(\C^2,0)$ and that  $\phi $ is analytically  isomorphic to  the map from $(\C^2,0)$  to $(\C^2,0)$ defined by $(x,y)  \mapsto (x,y^n).$
 
 2) The discriminant of $\phi $ is a normal crossing.  By definition $(X,0)$  is  then  a quasi-ordinary singularity  and   its  link is a lens space.    We prove that the minimal resolution of $(X,0)$ is a bamboo of rational curves (Proposition 5.3.1).

 Section 6 is an example of Hirzebruch-Jung's resolution.

 \subsection{ Conventions}

The boundary of a topological manifold $W$ will be denoted by $b(W).$

A {\bf disc}  (resp. an {\bf  open disc}) will always be an oriented   topological  manifold  orientation preserving homeomorphic  to $\{ z\in \C ,\vert z \vert \leq 1\}$ (resp. to $\{ z\in \C ,\vert z \vert < 1\}$).

A {\bf circle } will always be an oriented     topological  manifold orientation preserving homeomorphic  to $ S  = \{ z\in \C ,\vert z \vert =1\}.$ Moreover,  for $0<\alpha$,  we use the following notation:  $ D_{\alpha}  = \{ z\in \C ,\vert z \vert \leq \alpha \}$,  and $S_{\alpha }=b(D_{\alpha})$.

\bigskip
{\bf Acknowledgments: } 
 I thank Claude Weber for  useful discussions and  for making many  comments about the redaction of this text.

\newpage

\section{The topology of plumbings}

In this Section $(X,0)$  is  a  {\bf normal} complex surface germ.

The name ``plumbing" was introduced by David Mumford in \cite{mu} (1961). There, he showed  that the topology of a resolution of a normal singularity of a complex surface can be described as  a  "plumbing". 

 In  \cite{hi} (1953),  Hirzebruch  constructed  good  resolutions of normal singularities.
Let $ \rho  :  (Y, E_{Y}) \to  (X, 0)$ be a good resolution of the  normal germ of surface $(X,0).$ Each irreducible component $E_i$ of the exceptional divisor is equipped  with its normal complex fiber bundle. With their complex structure the fibers have dimension $1$. So, a  regular  compact tubular   neighbourhood $N(E_i)$  of $E_i$ in $Y$, is a  disc bundle. As $E_i$ is a smooth  compact complex curve, $E_i$ is an oriented  differential compact surface with an empty boundary. Then,  the  isomorphic class, as differential bundle,  of the disc bundle $N(E_i)$ is given by  the genus  $g_i$ of  $E_i$ and its self-intersection $e_i$ in $Y.$ The complex structure gives an orientation on $Y$ and  on  $E_i$, these orientations induce  an orientation on    $N(E_i)$  and on the fibers of the  disc bundle over $E_i.$

\begin{remark}
 By definition $(X,0)$ is a sufficiently  small compact representative of the given normal surface germ.   Let $k$ be the number of irreducible components of $E_Y$,    $ M(Y)= \cup _{1 \leq i \leq k}  N(E_i)$ is a compact neighborhood of $E_Y$. There exist  a retraction by deformation  $R: Y \to M(Y)$ which induces an homeomorphism from the boundary of Y,  $b(Y)= \rho ^{-1}(L_X)$, to the boundary $b(M(Y))$. So, the boundary of $M(Y)$ is the link of $(X,0).$ 
\end{remark}

\begin{definition}
Let $N(E_i),i=1,2,$  be two oriented  disc bundles on oriented smooth compact differentiable  surfaces,  with empty boundary, $E_i,i=1,2,$  and  let $p_i\in E_i$. {\bf The plumbing of  $N(E_1)$ and  $N(E_2)$ at $p_1$ and $ p_2$ } is equal to the  quotient of the disjoint union of  $N(E_1)$ and $N(E_2)$ by the following  equivalence relation.

Let $D_i$ be a small disc neighbourhood of $p_i$ in $E_i$,  and  $D_i \times \Delta_i$ be a trivialization of $N(E_i)$ over $D_i$, $i=1,2.$  Let $f : D_1 \to \Delta _2$ and $g: \Delta _1 \to D_2 $ be two orientation preserving diffeomorphisms such that $f(p_1)=0$ and $g(0)=p_2.$

 For all  $(v_1,u_1)\in D_1 \times \Delta _1$, the equivalence relation is $ (v_1,u_1) \sim (g(u_1),f(v_1))$. 
 
 \end{definition}
 
 \begin{remark}
 
 The diffeomorphism class of the plumbing of  $N(E_1)$ and  $N(E_2)$ at $(p_1, p_2)$ doesn't  depends upon  the choices  of the trivialisations  neither the choices of $f$ and $g.$ Moreover, in the plumbing  of  $N(E_1)$ and  $N(E_2)$ at $p_1$ and $ p_2$:
  \begin{enumerate}

 \item  The image of $E_1$ intersects  the image of $E_2$ at   the point  $p_{12}$ which is the class, in the quotient,    of   $(p_1 \times 0) \sim  (p_2 \times 0).$ 
 
 \item The plumbing is a gluing of 
 $N(E_1)$ and  $N(E_2)$  around the chosen neighhbourhoods of  $(p_1 \times 0)$ and $  (p_2 \times 0).$  
 \item  In the plumbing,  $D_1\times 0 \subset E_1$ is identified,   via $f$,  with the  fiber  $0\times  \Delta _2$   of the disc bundle $N(E_2)$ and the fiber $0\times \Delta _1$  of $N(E_1)$ is identified,   via $g$,  with $D_2 \times 0 \subset E_2.$  

 \end{enumerate}

 \end{remark}

\begin{definition} 

More generally we can perform the  plumbing of a family $N(E_i), i=1, \dots,  n,$ of oriented  disc bundles on oriented smooth compact differentiable  surfaces $E_i$ with empty boundary,  at  a finite number of  pairs of points  $(p_i, p_j) \in E_i \times E_j$.  Let $g_i$ be the genus of $E_i$ and $e_i$ be the self intersection of $E_i$ in $N(E_i)$.   The  vertices of  the  {\bf weighted  plumbing  graph} associated to  such a plumbing represent the basis $E_i, i=1, \dots ,n ,$ of the bundles. These  vertices  are weighted by $e_i$ when $g_i=0,$ and by $(e_i,g_i)$ when $0<g_i$.  Each edge which relates $(i)$ to $(j)$, represents an intersection  point between the image of $E_i$ and $E_j$ in the  plumbing. 
 
In the boundary of the plumbing of the family $N(E_i), i=1, \dots,  n$,  the intersections $b(N(E_i)) \cap b(N(E_j)$ are a union of disjoint tori which is the {\bf family of plumbing tori} of the plumbing.
 
\end{definition}

An  oriented  disc bundle  $N(E)$ on  a  differential  compact surface   $E$ of genus $g$ and empty boundary   is determined  as differentiable bundle by $g$ and  by the self-intersection of $E$ in $N(E)$. 
If two plumbings have the same  weighted  plumbing graph, there exists a diffeomorphism between the two plumbings such that its restriction on the corresponding disc bundles is an isomorphism of differentiable  disc bundles.

\begin{proposition} 
Let $ \rho  :  (Y, E_{Y}) \to  (X, 0)$ be a good resolution of the  normal germ of surface $(X,0).$ Then a regular neighbourhood, in $Y,$  of the exceptional divisor divisor $E_{Y}$,  is diffeomorphic to a plumbing of the disc bundles  $N(E_i)$. The plumbings are performed around the double points $p_{ij}= E_i \cap E_j.$The associated  weighted plumbing graph is equal to the weighted  dual graph  $G_w(Y)$ of $\rho$. To each point $p_{ij} \in (E_i \cap E_j )$ we associate a torus $T(p_{ij}) \subset  (b(N(E_i)) \cap b(N(E_j))).$
\end{proposition}

 {\it  Proof:}  We choose  trivializations of the disc bundles $N(E_i)$ and $N(E_j)$ in a small closed neighborhood $V$ of  $p_{ij}.$ First,  we centered the trivializations at $(0,0)=p_{ij}$,  we parametrized $V $ as disc bundle  
 
 1) over  $E_i$ by $V_i=\{ (v_i,u_i)\in D_i \times \Delta _i\} $, where $D_i \times 0$ is a disc neighborhood of   $(0,0)=p_{ij}$ in $E_i$ and $v_i\times \Delta _i$ is  the normal disc fiber at $v_i\in D_i$. 
 
  2) over  $E_j$ by $V_j=\{ (v_j,u_j)\in D_j \times \Delta _j \} $, where $D_j\times 0$ is a disc neighborhood of   $(0,0)=p_{ij}$ in $E_j$ and $v_j \times \Delta _j $ is  the normal disc fiber at $v_j \in D_j.$

  As $E_Y$ is a normal crossing divisor, we can parametrize $V$ such that  $ E_Y \cap V =\{uv=0\}$
  where $v=v_i=u_j$ and $u=v_j=u_i$. These equalities provide the plumbing of $N(E_i)$ and $N(E_j)$ around $p_{ij}.$

    By construction , the  associated  weighted plumbing  graph is equal  to $G_w(Y).$
    
     \begin{definition}  Now we can say that the union of disc bundles  $ M(Y)= \cup _{1 \leq i \leq k}  N(E_i)$   is {\bf the plumbing}   associated to $ \rho  :  (Y, E_{Y}) \to  (X, 0)$. 
     \end{definition}

 With the above notations, in a neighborhood of $p_{ij}$, there is a unique connected component of the intersection  $(b(N(E_i)) \cap b(N(E_j)))$  which is parametrized by  the torus  $b(D_i) \times b(\Delta _i)$ which is glued point by point with  $b(D_j) \times b(\Delta _j)$. 
 
  \begin{definition}
  The image of $(b(D_i) \times b(\Delta _i)) \sim (b(D_j) \times b(\Delta _j))$ in the boundary of $M(Y)$ is the {\bf plumbing torus $T(p_{ij})$ associated to $p_{ij}$.}

\end{definition}

{\it End of proof.}

\newpage

\section{The topology of the normalization}

In this Section $(X,0)$ is a reduced complex surface germ which can have a 1-dimensional singular locus.

\subsection{ $L_X$ as singular  covering over $S^3 $}

We choose a general projection   $\pi :  (X,0) \to  (\C ^2, 0)$, we denote by $\Gamma $ the singular locus of $\pi$ (in particular $\Sigma \subset \Gamma$) and by $\Delta $ its discriminant ($ \Delta = \pi (\Gamma) ) .$   In fact it is sufficient to  choose new coordinates in  $\C^n$,   $(x,y, w_1,\dots ,w_{n-2}) \in \C ^n$,    such that       the restriction  on $(X,0)$   of the projection 
$$(x,y, w_1,\dots ,w_{n-2}) \mapsto  (x,y),$$ 
 denoted by $\pi$,      is finite. Moreover, we assume that     the hyperplanes  $ H_a = \{x=a\} $ where $\vert a \vert \leq \alpha << \epsilon, $ meet  transversally  the singular locus $\Gamma$ of $\pi.$ In particular, $H_0\cap \Gamma = \{0\}.$

{\bf Conventions and notations}

Let $D_{\alpha }\times D_{\beta} \in \C^2$ be a polydisc at the origin  in $\C^2$  where   $0<\alpha <\beta<\epsilon$  are chosen sufficiently  small such that: 

I)  ${\cal{B}} = B_{\epsilon}^{2n} \cap \pi ^{-1}(D_{\alpha }\times D_{\beta})$ is a good semi-analytic neighborhood of $(X,0)$ in the sense of A. Durfee \cite{du}. Then $(X\cap {\cal B},0)$ is homeomorphic to $(X,0).$ In this section $(X,0)$  is given by  $(X\cap {\cal B},0)$. The link  $L_X=X\cap  b({\cal{B}})$ is the link of $X$. The link of $\Gamma $ is the link $K_{\Gamma}=\Gamma \cap  b({\cal{B}}) $  embedded in $L_X$.

II) For a sufficiently small $\alpha $, the intersection $ K_{\Delta}=\Delta \cap ((S_{\alpha}\times D_{\beta}) \cup (D_{\alpha }\times S_{\beta}))\subset  (S_{\alpha}\times D_{\beta}) .$ In this Section,  we choose such a  $ K_{\Delta}$ to represent the link of $\Delta$ embedded in  the 3-sphere (with corners)
$((S_{\alpha}\times D_{\beta}) \cup (D_{\alpha }\times S_{\beta}) )$.

III)  Let $\delta _j, 1\leq j \leq r, $ 
be the $r$ branches of  the discriminant $\Delta$. Let $N( K_{\Delta} )$ be a tubular compact neighborhood  of $ K_{\Delta}$.
So,  $N( K_{\Delta} )$ is a disjoint union of $r$ solid tori. For  a sufficiently small  $N( K_{\Delta} )$, the union  $N(K_{\Gamma })$,  of the connected components of $L_X \cap \pi^{-1} (N( K_{\Delta} ))$ which contain a connected component of  $K_{\Gamma}$, constitutes a tubular compact neighbourhood  of $K_{\Gamma}$ in $L_X$.

Let us denote by $\mathring{N}( K_{\Delta} ) $ the interior of $N( K_{\Delta} )$.    {\bf The exterior  $M$ of the link $K_{\Delta}$}   is  define by:  $$  M=((S_{\alpha}\times D_{\beta}) \cup (D_{\alpha }\times S_{\beta}) )\setminus \mathring{N}(K_{\Delta}) .$$

IV) Let $\gamma $ be a branch of  the singular locus $\Gamma $ of $\pi $. So,   $\pi (\gamma ) = \delta  $ is a branch of $\Delta $. Let $N(K_{\delta  })$ ( resp. $N(K_{\gamma })$) be the connected  component of $N(K_{\Delta })$ (resp. of $N(K_{\Gamma })$) which contains the link $K_{\delta }$  (resp.  $K_{\gamma }$).

\begin{remark} The restriction $\pi_L :L_X \to ((S_{\alpha}\times D_{\beta}) \cup (D_{\alpha }\times S_{\beta}) )$ of $\pi$ on $L_X$ is a finite morphism, its   restriction  on $M$ is a finite regular covering.  If $\gamma $ is not a branch of the singular locus $\Sigma $ of $X$, $\pi_L$ restricted to $N(K_{\gamma })$  is a ramified covering with $K_{\gamma }$ as  ramification locus.  If $\gamma $ is a branch of $\Sigma $,   $N(K_{\gamma })$  is  a singular pinched solid torus  as defined  in 3.4 and $\pi_L$ restricted to  $N(K_{\gamma })$  is   singular  all along $K_{\gamma }$.

\end{remark}

\subsection{\bf  Waldhausen  graph manifolds and plumbing graphs}

\begin{definition}
A Seifert fibration on an oriented, compact 3-manifold is an oriented foliation  by circles such that every leaves  has a tubular neighbourhood (which is a solid torus) saturated by leaves. A Seifert 3-manifold is an oriented, compact 3-manifold equipped with a Seifert fibration.
\end{definition}

\begin{remark} 

 \begin{enumerate}
\item  A Seifert 3-maniflod  $M$ can have a non-empty boundary. As this boundary is equipped  with a foliation by circles,  if $B(M)$ is non-empty it is a disjoint union of tori.

\item  Let $D$ be a disc and $r$ be a rotation of angle $2\pi q/p$ where $(q,p)$ are two positive integers prime to each other and $0 < q/p<1$. Let $T_r$  be the solid torus  equipped with a Seifert foliation given by the trajectories    of $r $ in the following mapping torus:  $$T_r=D\times \lbrack 0,1\rbrack / (z,1)\sim (r(z),0) .$$  
In particular, $l_0=  (0\times \lbrack 0,1\rbrack) / (0,1)\sim (0,0) $ is a core of $T_r.$ The other leaves are $(q,p)$-torus knots in $T_r$.  Let $T_0$ be $D\times S$ equipped with the trivial fibration by circles $l(z)=\{z\} \times S, \ z\in D$.  A solid torus  $T(l)$   which is a tubular neighbourhood of a leave $l$ of a Seifert 3-manifold $M$ is either  

1) orientation and foliation preserving homeomorphic to  $T_0$. In this case,  $l$ is a regular Seifert  leave.

2) or,  is orientation and foliation preserving homeomorphic to $T_r$. In this case, $l$ is an exceptional leave of $M$.

\item  The compacity of $M$ implies that there exists  only a finite number of exceptional leaves.

 \end{enumerate}

\end{remark}

\begin{definition}

Let $M$ be  an oriented and  compact 3-manifold. The manifold  $M$ is a {\bf  Waldhausen graph manifold}  if there exists  a finite family $\cal{T} $,  of disjoint tori  embedded in $M$,  such that if $M_i, i=1,\dots , m,$  is the family of the closure of the connected components of $M\setminus \cal{T}$, then  $M_i$ is a Seifert manifold for all $i, \ 1\leq i \leq m$. We assume that it  gives us   a  finite decomposition $M = \cup _{1\leq i \leq m} M_i$ into a union of compact connected  Seifert manifolds  which satisfies  the following properties:

\begin{enumerate}
\item  For each $M_i, i=1,\dots ,m$,  the boundary of $M_i$ is in $\cal{T} $ i.e. $b(M_i)\subset \cal{T} $.
\item If $i \neq j$ we have the inclusions  $(M_i \cap M_j) \subset \cal{T} $. 
\item The  intersection   $(M_i \cap M_j) $,  between  two  Seifert manifolds of the given decomposition,  
 is either empty or  equal to the  union of the common boundary components of $M_i$ and $M_j$.
\end{enumerate}

 Such a decomposition  $M = \cup _{1\leq i \leq m} M_i$,  is   {\bf the Waldhausen decomposition of $M$,  associated to  the family of tori $\cal{T} $}.
  
\end{definition}

\begin{remark}

One can easily deduce from 2.4,  that the  family of the  plumbing tori    gives a decomposition of the  boundary of a plumbing as a union of Seifert manifolds. So, the boundary of a plumbing is a Waldhausen  graph manifold.

In \cite{ne}, W.  Neumann shows how  to construct a plumbing from a  given Waldhausen decomposition of a 3-dimensional oriented compact manifold. 
\end{remark}

  As in Section 3.1, we consider  the exterior  $  M=((S_{\alpha}\times D_{\beta}) \cup (D_{\alpha }\times S_{\beta}) )\setminus \mathring{N}(K_{\Delta}) $ of 
the link $ K_{\Delta}$.  The  following proposition is well known  (for example see \cite{e-n}, \cite{m-w}). Moreover,   a detailed description of $M$,  as included in  the boundary of the plumbing graph given by the minimal resolution of $\Delta$,  is given in \cite{l-m-w}, p. 147-150.

\begin{proposition} The exterior  $M$ of the link of a plane curve germ $\Delta$ is a Waldhausen graph manifold.  The minimal   Waldhausen  decomposition of $M$ can be extended to a Waldhausen decomposition of the sphere $((S_{\alpha}\times D_{\beta}) \cup (D_{\alpha }\times S_{\beta}) )$  in which the    connected components of $K_{\Delta}$ are  Seifert leaves. Moreover, with such a Waldhausen decomposition,   the solid tori connected components of $N(K_{\Delta})$ are saturated by  Seifert leaves which are oriented circles transverse to $(a\times D_{\beta}), a \in S_{\alpha}$. The   cores $K_{\Delta}$  of $N(K_{\Delta})$  are  a union of these Seifert leaves.

\end{proposition}

\subsection {\bf The topology of $L_X$ when $L_X$ is a topological manifold }

 If $(X,0)$ is not normal, let    $\nu_L : L_{\bar X} \to  L_X$  be the normalization of $(X,0)$ restricted to the link of $(\bar {X},p)$ (if $(X,0)$ is normal $\nu_L$ is the identity). 

\begin{remark}  The link of a normal  complex surface germ is a  Waldhausen  graph manifold. Indeed, the composition morphism $\pi_L\circ \nu_L$ is a ramified covering with the link $K_{\Delta}$ as set of ramification values:
$$(\pi_L\circ \nu_L): L_{\bar X} \to ((S_{\alpha}\times D_{\beta}) \cup (D_{\alpha }\times S_{\beta}) ).$$
So, we can take  inverse image, by $(\pi_L\circ \nu_L)^{-1}$,  of the tori and of the Seifert leaves of a Waldhausen decomposition of   $((S_{\alpha}\times D_{\beta}) \cup (D_{\alpha }\times S_{\beta}) )$   in  which  $K_{\Delta}$ is a union of Seifert leaves, to obtain  a Waldhausen decomposition of $ L_{\bar X}$. Then, the   plumbing calculus \cite{ne}, describes $L_{\bar X}$ as the boundary of a plumbing without the help of a good  resolution of $(\bar X,p)$.

\end{remark}

If the singular locus $(\Sigma,0)$ of $(X,0)$ is one dimensional, let $(\sigma,0) $ be a branch of $(\Sigma,0) $  and  $s$ be a point of the intersection  $\sigma \cap \{x=a\} $. Let $\delta = \pi (\sigma) $ be the branch of  the discriminant $\Delta $ which is the image  of $ \sigma $ by the morphism $\pi $. Then,  $\pi _L  ( s)=(a,y)\in   (S_{\alpha}\times D_{\beta}) $. Let $N(K_{\delta })$ be a solid torus regular neighbourhood of $K_{\delta }$ in $(S_{\alpha}\times D_{\beta})$ and  let $N(K_{\sigma  })$ be the connected component of $(\pi _L)^{-1}( N(K_{\delta })) $ which contains  $s$ (and $K_{\sigma }$). 

\begin{definition} 
 \begin{enumerate}

 \item  Let $(C,s)$ be the germ of curve which is the connected component of  $ N(K_{\sigma }) \cap \{x=a\} $ which contains $s$.    By definition $(C,s)$ is the {\bf hyperplane section germ} of $\sigma $ at $s$. For a sufficiently small $\alpha =\vert a \vert $, $(C,s)$  is reduced and  its topological type does not depend upon  the choice of $s$. In particular, $k(\sigma )$,  the number  of the  irreducible components of  $(C,s)$,  is well defined.

If $k(\sigma )=1$, $\sigma $ is a {\bf  branch} of $\Sigma $ {\bf with irreducible hyperplane sections}. Let $\Sigma =\Sigma _1 \cup \Sigma _+$ where $\Sigma_1$ is the union of the branches of $\Sigma $ with irreducible hyperplane sections and  where $\Sigma _+$  is the union of the branches of $\Sigma $ with reducible  hyperplane sections. 

 \item Let $D_i, 1\leq  i\leq k$ be $k$   oriented discs centered at $0_i\in D_i$. A  {\bf $k$-pinched disc $k(D)$ } is a topological space orientation preserving  homeomorphic to the quotient of the disjoint union of the $k$ discs by the identification of their centers in a unique   point $\tilde{0}$ i.e. $0_i \sim 0_j$ for all $i$ and $j$ where $ 1\leq  i\leq k,1\leq  j\leq k$. {\bf The center of $k(D)$} is  the  equivalence class $\tilde{0}$ of the centers $0_i,1\leq  i\leq k$.
 
  \item If  $\ h: k(D) \to k(D)' $ is an  homeomorphism between two $k$-pinched discs with $k>1$, $h(\tilde{0})$ is obviously  the center of $k(D)'$. We say that {\bf $h$ is orientation preserving} if $h$ preserves the orientation of the punctured $k$-pinched discs  $(k(D)\setminus \{ \tilde{0} \})$ and $(k(D)'\setminus \{ \tilde{0} \})$.
 
  \end{enumerate}
 
\end{definition}

\begin{lemma} Let $(C,s)$ be the germ of curve which is the connected component of  $ N(K_{\sigma }) \cap \{x=a\} $ which contains $s$. Then, $C$ is a $k(\sigma )$-pinched disc centered at $s$ and 
$N(K_{\sigma })$ is the mapping torus of $C$ by an orientation preserving homeomorphism $h$ which fixes  the point $s$. 

\end{lemma} 
{\it Proof:}  As $(C,s)$ is a germ of curve with $k(\sigma )$ branches, up to homeomorphism $( C,s)$ is a $k(\sigma )$-pinched disc with $s=\tilde{0}$.

We can saturate  the solid torus $ N(K_{\delta })=\pi ( N(K_{\sigma }))$ with oriented circles such that $K_{\delta }$ is one of these circles and such that  the first return  homeomorphism defined  by these circles on the disc  $\pi ( C)$ is   the  identity. Let $\gamma $ be one circle of the chosen saturation of $ N(K_{\delta })$, $\pi^{-1}(\gamma)\cap  N(K_{\sigma })$ is a disjoint union of oriented circles  because $\pi $ restricted to $N(K_{\sigma })\setminus K_{\sigma }$ is a regular covering and $(\pi^{-1}(K_{\delta})\cap  N(K_{\sigma }))=K_{\sigma}$. So, $N(K_{\sigma })$ is equipped with a saturation by oriented circles. The first return map on $C$ along the so constructed circles is an orientation preserving homeomorphism $h$ such that $h(s)=s$ because $K_{\sigma }$ is one of the given circles.

{\it End of proof.}

\begin{lemma} As  above, let $(C,s)$ be the hyperplane section germ at $s\in \sigma \cap \{x=a\} $.
Let  $\bar {\sigma}_j,1\leq j \leq n,$ be the $n$ irreducible components  of $\nu _L^{-1}(\sigma)$ and let $d_j$ be the order of   $\nu_L$ restricted to $\bar {\sigma}_j $. Then, we have $$k(\sigma)=d_1+\dots +d_j+\dots +d_n.$$
\end{lemma}
{\it Proof:}  The normalization $\nu$ restricted to $\bar {X} \setminus \bar {\Sigma }$, where $\bar {\Sigma }= \pi ^{-1}(\Sigma)$, is an isomorphism. The number $n$ of  the irreducible components  of $\nu _L^{-1}(\sigma)$ is equal to the number of the connected components of $\nu _L^{-1}(N(K_{\sigma}))$. So, $n$ is the number of the connected components of  the boundaries   $b(\nu _L^{-1}(N(K_{\sigma})))$ which is equal to the number of the connected components of    $b(N(K_{\sigma}))$. Let $\tau_j, 1\leq j \leq n,$ be the  $n$ disjoint tori which are the  boundary  of   $N(K_{\sigma})$. The order $d_j$ of $\nu$ restricted to $\bar {\sigma}_j$ is equal to the number of points of  $\nu _L^{-1} (s) \cap (\bar {\sigma}_j )$. 
 
 Let $(\gamma _j,s)$ be an irreducible component of $(C,s)$ such that $m_j= b(\gamma _j) \subset \tau_j $. The normalization   
  $\nu $ restricted to $ (\nu_L^{-1}(\gamma_j \setminus \{s\} ))$ is an isomorphism over the punctured disc $(\gamma_j \setminus \{s\} )$.  So, the intersection      $\nu _L^{-1} (\gamma_j ) \cap  \bar {\sigma}_j$ is a unique point $p_j$. As $(\bar {X}, p)$ is normal, $p_j$ is  a smooth point of $(\bar {X},p)$ and then, $\nu _L^{-1} (\gamma_j ) $ is  irreducible  and it  is the only irreducible component of  $\nu _L^{-1} (C ) $ at $p_j$. By symmetry there is exactly one irreducible component of  $\nu _L^{-1} (C ) $ at every point  of $\nu _L^{-1} (s) \cap (\bar {\sigma}_j )$.

 So, $d_j$ is the number of the meridian circles  of the solid torus $N(K_{\bar {\sigma}_j})$ obtained by  the following  intersection  $ ( \nu _L^{-1} (C )) \cap (\nu _L^{-1}( \tau _j))$. But $\nu $ restricted to 
$(\nu _L^{-1}( \tau _j))$ is an isomorphism and $d_j$ is also the number of connected components of $C \cap \tau_j$. So,  $d_1+\dots +d_j+\dots +d_n,$ is   equal to the number of connected components of $b (  C)=C\cap b(N(K_{\sigma}))$ which is the number of irreducible components of $(C,s)$.   

{\it End of proof.}

\begin{remark} A well-known result of analytic geometry could be roughly  stated as follows: ``The normalization separates the irreducible components". Here, $(X,0)$ has $k(\sigma)$ irreducible components around  $s\in  \sigma$.    Using only  basic topology, Lemma 3.3.4 proves that $ ( \nu _L^{-1} (s ))$ has $k(\sigma)=d_1+\dots +d_j+\dots +d_n$ distinct points  and that there is exactly one irreducible component of  $\nu _L^{-1} (C ) $ at every point  of $\nu _L^{-1} (s)$. This gives a topological proof that the normalization $\nu$ separates the  $k( \sigma)$ irreducible components of $(C,s)$ around  
$s\in  \sigma$.

\end{remark}

\begin{proposition} The  following three statements are equivalent:
\begin{enumerate} 
\item $L_X$ is a topological manifold equipped with a   Waldhausen graph  manifold structure.
\item The normalization $\nu: (\bar{X} ,p) \to (X,0) $ is a homeomorphism.
\item All the branches of $\Sigma $ have irreducible hyperplane sections.
 
 \end{enumerate}
\end{proposition}

{\it Proof:}  The normalization $\nu$ restricted to $\bar {X} \setminus \bar {\Sigma }$, where $\bar {\Sigma }= \pi ^{-1}(\Sigma)$, is an isomorphism. The normalization is a homeomorphism if and only if $\nu$  restricted to  $\bar {\Sigma }= \pi ^{-1}(\Sigma)$ is a bijection. This  is the case if and only if  we have   $1=d_1+\dots +d_j+\dots +d_n \ $  for all the branches $\sigma $ of $\Sigma $. But, by Lemma 3.3.4,    $k(\sigma )=d_1+\dots +d_j+\dots +d_n.$ This proves the equivalence of the statements 2 and 3.
 
Let $(C,s)$ be the hyperplane section germ at $s\in \sigma \cap \{x=a\} $.  If $L_X$ is a topological manifold, it is a topological manifold at $s$  and $k(\sigma )=1$ for all branches $\sigma$ of $\Sigma$.  
 If all the branches of $\Sigma $ have irreducible hyperplane sections, we already know that the normalization $\nu: (\bar{X} ,p) \to (X,0) $ is a homeomorphism. Then, the  restriction $\nu_L$ of $\nu$ on $L_{\bar X}$ is also  a homeomorphism.  By Remark 3.3.1,  $L_{\bar X}$ is a Waldhausen graph manifold. In particular, we can equip  $L_X$ with the Waldhausen graph manifold  structure carried by $\nu_L$. This proves the equivalence of the statements 1 and 3.

{\it End of proof.}

\subsection{ Singular $L_X$, curlings and identifications}

In 3.3 (Definition 3.3.2), we have considered the union  $\Sigma_+$ of the branches of the singular locus $\Sigma $ of $(X,0)$ which have reducible hyperplane sections. We consider a tubular neighbourhood    $N_+=\cup _{\sigma \subset   \Sigma_+} N(K_{\sigma})$ of the link  $K_{\Sigma _+}$ of $\Sigma_+$  in $L_X$. As in the proof of Proposition 3.3.4,  the exterior $M_1= L_X \setminus  \mathring {N}_+$,  of $K_{\Sigma _+}$  in $L_X$, is a topological manifold because $\nu $ restricted to 
$\nu ^{-1} (M_1)$ is a homeomorphism. From now on  $\sigma$ is a branch of $\Sigma _+$. By  Proposition 3.3.4, $L_X$ is topologically  singular at every point of  $K_{\sigma}$. In this section, we show that $N(K_{\sigma})$ is a singular pinched  solid torus. In Lemma 3.3.3, it is shown that $N(K_{\sigma})$ is the mapping torus of a $k(\sigma )$-pinched disc  by an orientation preserving homeomorphism. But, the homeomorphism class of the  mapping torus  of a homeomorphism $h$ depends only on the isotopy class of $h$. Moreover the isotopy class of an orientation preserving homeomorphism  $h$ of a $k$-pinched disc depends only on the permutation induced by $h$ on the $k$ discs. In particular,  if  $h: D\to D$ is an orientation preserving  homeomorphism of a disc $D$   the associated  mapping torus $$T(D, h)= \lbrack 0,1 \rbrack \times D / (1,x) \sim (0,h(x))$$ is homeomorphic to the standard torus $S\times D$.

\begin{definition}
 \begin{enumerate}
 \item   Let  $k(D)$  be the $k$-pinched disc quotient by identification  of their centrum of $k$ oriented  and ordered discs $D_i, 1\leq  i\leq k$. Let $c=c_1 \circ c_2  \circ  \dots  \circ c_n$  be a permutation of  the indices $\{1, \dots ,k\}$  given as the composition of $n$ disjoint cycles $c_j, 1\leq j \leq n,$  where $c_j$ is a cycle of order $d_j$.  Let $\tilde{h}_c$ be an orientation preserving homeomorphism of the disjoint union of  $D_i, 1\leq  i\leq k$ such that $\tilde{h}_c (D_i)=D_{c(i)}$ and  $ \tilde{h}_c(0_i)=0_{c(i)}$.  Then,  $\tilde{h_c}$ induces  an  orientation preserving homeomorphism  $h_c$  on  $k(D)$.  By construction we have  $h_c(\tilde{0})= \tilde{0}$.    A  {\bf singular pinched  solid torus  associated to the permutation $c$} is an topological space orientation preserving  homeomorphic to the  mapping torus  $T(k(D),c)$ of $h_c$:
 $$ T(k(D),c)=  \lbrack 0,1 \rbrack  \times k(D) / (1,x) \sim ( 0,h_c(x))$$
 
  The {\bf core } of $T(k(D),c)$ is the  oriented circle    $ l_0=   \lbrack  0,1 \rbrack \times  \tilde{0}  / (1,\tilde{0}) \sim   (0,\tilde{0})$.  A  homeomorphism between two singular pinched  solid tori is orientation preserving if it preserves the orientation of  $k(D)\setminus \{ \tilde{0} \}$ and the orientation of  the trajectories of $h_c$  in  its mapping torus $ T(k(D),c)$.
  
 \item  A  {\bf $d$-curling  $\cal {C}_d$} is  a topological space  homeomorphic to    the following quotient of a solid torus $S\times D:$
 $$\cal {C}_d=S\times D /(u, 0)\sim (u',0) \Leftrightarrow u^d=u'^d.$$
 Let $q:(S\times D)\to \cal {C}_d$  be the associated quotient morphism. By definition,   $l_0=q(S\times \{0\})$ is the {\bf core of $\cal {C}_d$}.
  \end{enumerate}

\end{definition}

\begin{example}
Let $X= \{(x,y,z)\in \C ^3 \ where \  z^d-xy^d=0\}$. The normalization of $(X,0)$ is smooth i.e. $\nu : (\C^2,0) \to (X,0)$ is given by $(u,v) \mapsto  (u^d,v,uv)$. Let $T= \{(u,v)\in  (S \times D ) \subset \C^2 \}$.   Let $\pi_x : \nu (T) \to S $ be the projection  
$(x,y,z) \mapsto   x$ restricted to $\nu (T)$. Here the singular locus of $(X,0)$ is the line $\sigma =(x,0,0), x\in \C$.   We have   $N(K_{\sigma}) = L_X \cap (\pi_x^{-1}(S)) =\nu (T)$ as  a tubular neighbourhood of $K_{\sigma}$.  Let  $q: T \to \cal {C}_d$ be the quotient morphism defined above.
There exists  a well defined homeomorphism  $f: \cal {C}_d \to N(K_{\sigma})$ which satisfies 
$f(q(u,v))=(u^d,v,uv)$. So,  $N(K_{\sigma})$ is a    d-curling  and  $K_{\sigma}$ is its core. Moreover, $f$ restricted to the  core $l_0$  of $ \cal {C}_d $  is a homeomorphism onto  $K_{\sigma}$.

\end{example}

\begin{figure}[ht]

\begin{tikzpicture} 

\begin{scope}[xshift=2cm]
\begin{scope}[yshift=10cm]

 \draw[line width=1pt, lightgray  ] 
(5,3.5)  .. controls (5.5,3.5)  and  (5.1,2.8) .. (5.1,2.6) ;  
 \draw[line width=1pt, lightgray  ] 
(5.8,1.2) .. controls (5.8,1.7) and (5.1,2.3)  .. (5.1,2.6) ;
 \draw[line width=2pt] 
(5.19,2.95)  .. controls (5.08,2.65)  and  (5.08,2.6) .. (5.23,2.3) ;

 \draw[line width=1pt, lightgray   ] 
(4.7,2.3) .. controls (4.7,2.6) and (4.5,3.5)  .. (5,3.5) ;

\draw[line width=1pt, lightgray  ] 
(4.7,2.3) .. controls (4.7,2) and (4.2,1.6)  .. (4.2,1.1) ;

 \draw[line width=2pt] 
(4.67,2.7)  .. controls (4.73,2.5)  and  (4.75,2.2) .. (4.52,1.9) ;

4.7,2.3

\draw[line width=1pt, lightgray  ] 
(5,0) .. controls (5.7,0.1) and (5.8,0.7)  .. (5.8,1.2) ;

 \draw[line width=1pt , lightgray ]
(4.2,1.1) .. controls (4.2,1) and (4,0)  .. (5,0) ;



\draw[line width=0.5pt] 
(1,1.5) .. controls (2,2.5) and (4,3.4)  .. (5,3.51) ;

\draw[line width=0.5pt] 
(1,1.5) .. controls (2,2) and (4,2.3)  .. (4.7,2.3) ;

\draw[line width=0.5pt] 
(4.65,2.62) .. controls (4.9,2.65) and (5,2.65)  .. (5.1,2.64) ;

\draw[dotted, line width=0.6pt] 
(1,1.5) .. controls (2,2.1) and (4,2.64)  .. (4.65,2.62) ;

\draw[line width=0.5pt] 
(1,1.5) .. controls (2,1.5) and (4.3,0.3)  .. (4.5,0.1) ;

\draw[fill=white] (4.7,2.3)circle(1.5pt);
\draw[fill=white] (5.11,2.64)circle(1.5pt);

 \node at(0.8,1.5)[ ]{$\bar{p}$};

   \node at(6,2.3)[right]{$\overline{L_X}$};
 \draw[line width=0.3pt, lightgray]   (6,2.3)--(5.25,2.85);
  \draw[line width=0.3pt, lightgray]   (6,2.25)--(4.7,2);
       
       \node at(2.6,3.5)[left]{$\overline{K}$};
        \draw[line width=0.3pt, lightgray]   (2.5,3.35)--(4.6,2.4);
         \draw[line width=0.3pt, lightgray]   (2.5,3.4)--(5,2.7);
        
  \node at(3.5,1.2)[left]{$\overline{X}$};
  
  \end{scope}

 \begin{scope}[yshift=4.5cm]
  \node at(3.1,4.8)[right]{$\nu $};
  \draw[line width=0.5pt]   (3,5.2)--(3,4.5);
   \draw[line width=0.5pt]   (3.1,4.7)--(3,4.5);
     \draw[line width=0.5pt]   (2.9,4.7)--(3,4.5);
     \end{scope}

  \begin{scope}[yshift=5cm]

 \draw[line width=1pt, lightgray  ]
(4.2,1.1) .. controls (4.2,1) and (4,0)  .. (5,0) ;

\draw[line width=1pt, lightgray ] 
(5,2.5) .. controls (4.8,2.2) and (4.2,1.6)  .. (4.2,1.1) ;

 \draw[line width=1pt, lightgray    ] 
(5,0) .. controls (6,0.1) and (6,1.8)  .. (5,2.5) ;

 \draw[line width=1pt, lightgray   ] 
(5,2.5) .. controls (4.5,2.9) and (4.6,3.5)  .. (5,3.5) ;
\draw[line width=2pt] 
(5,2.5) .. controls (4.9,2.6) and (4.8,2.7)  .. (4.75,2.8) ;
\draw[line width=2pt] 
(5,2.5) .. controls (5.1,2.42) and (5.05,2.47)  .. (5.29,2.23) ;
 
 \draw[line width=1pt, lightgray    ] 
(5,3.5)  .. controls (5.2,3.5)  and  (5.4,3.2) .. (5,2.5) ;

\draw[line width=2pt] (5.17,2.85)  .. controls (5.05,2.6)  and  (5.05,2.6) .. (5,2.5) ;
\draw[line width=2pt] (5,2.5) .. controls (4.95,2.4) and (4.85,2.3)  .. (4.72,2.15) ;

\draw[line width=0.5pt] 
(1,1.5) .. controls (2,2.5) and (4,3.4)  .. (5,3.5) ;

\draw[line width=0.5pt] 
(1,1.5) .. controls (2,2) and (4,2.5)  .. (5,2.5) ;

\draw[line width=0.5pt] 
(1,1.5) .. controls (2,1.5) and (4.3,0.3)  .. (4.5,0.1) ;

  \draw[fill=white] (5,2.5)circle(1.5pt);

 \node at(0.8,1.5)[ ]{$0$};

   \node at(6,2.3)[right]{$L_X$};
 \draw[line width=0.3pt, lightgray]   (6,2.3)--(5.13,2.65);

       \node at(2.6,3.5)[left]{$K$};
        \draw[line width=0.3pt, lightgray]   (2.5,3.4)--(4.85,2.54);
  
   \node at(0.8,1.5)[ ]{$0$};
  \node at(3.6,1.8)[ ]{$X$};
   \node at(2,2.5)[ ]{$\gamma_1$};
    \node at(2.8,2.3)[ ]{$\gamma_2$};
      \node at(3,1.1)[ ]{$\gamma_3$};

   \node at(8,2.5)[right ]{$\Sigma=\gamma_2 $};
    \node at(8,2)[ right]{$\Gamma = \bigcup_{i=1}^3 \gamma_i$  };

\end{scope}
 
\begin{scope}[yshift=-0.5cm]
  \node at(3.1,4.8)[right]{$\pi $};
  \draw[line width=0.5pt]   (3,5.2)--(3,4.5);
   \draw[line width=0.5pt]   (3.1,4.7)--(3,4.5);
     \draw[line width=0.5pt]   (2.9,4.7)--(3,4.5);
     \end{scope}
     
 \begin{scope}[xshift=2.5cm]
 
 \node at(0.8,1.5)[ ]{$0$};

\draw[line width=0.5pt] 
(1,1.5) .. controls (2,2.5) and (4,3.4)  .. (5,3.5) ;

\draw[line width=0.5pt] 
(1,1.5) .. controls (2,2) and (4,2.5)  .. (5,2.5) ;

 \draw[line width=0.5pt] 
(1,1.5) .. controls (2,1.5) and (4.3,0.3)  .. (4.5,0.1) ;

\draw[line width=0.3pt] 
(-0.5,3.5)--(-0.5,-0.5)--(2.5,-0.5)--(2.5,3.5)--(-0.5,3.5);

\draw[line width=2pt] 
 (2.5, 2.4)--(2.5,2.8);
 
 \draw[line width=2pt] 
 (2.5, 2.2)--(2.5,1.9);
 
   \draw[line width=2pt] 
 (2.5, 1.3)--(2.5,0.8);

   \node at(3.4,3.25)[ ]{$\delta_1$};
 \node at(3.4,2.5)[ ]{$\delta_2$};
 \node at(3.4,0.95)[ ]{$\delta_3$};
      
        \node at(2.25,0.9)[ ]{$V$};
         \node at(-0.5,-0.2)[right ]{$ D^2_{\alpha} \times  D^2_{\beta}$};
          
        \node at(5.5,1.5)[right ]{$\Delta=\bigcup_{i=1}^3\delta_i$};
        
  \end{scope}

\end{scope}

 \end{tikzpicture} 
  \caption{Schematic picture of $\pi $ and  $ \nu$ when  there is a 2-curling on $\Sigma = \gamma_2$}
  \label{fig:nor}
  \end{figure}
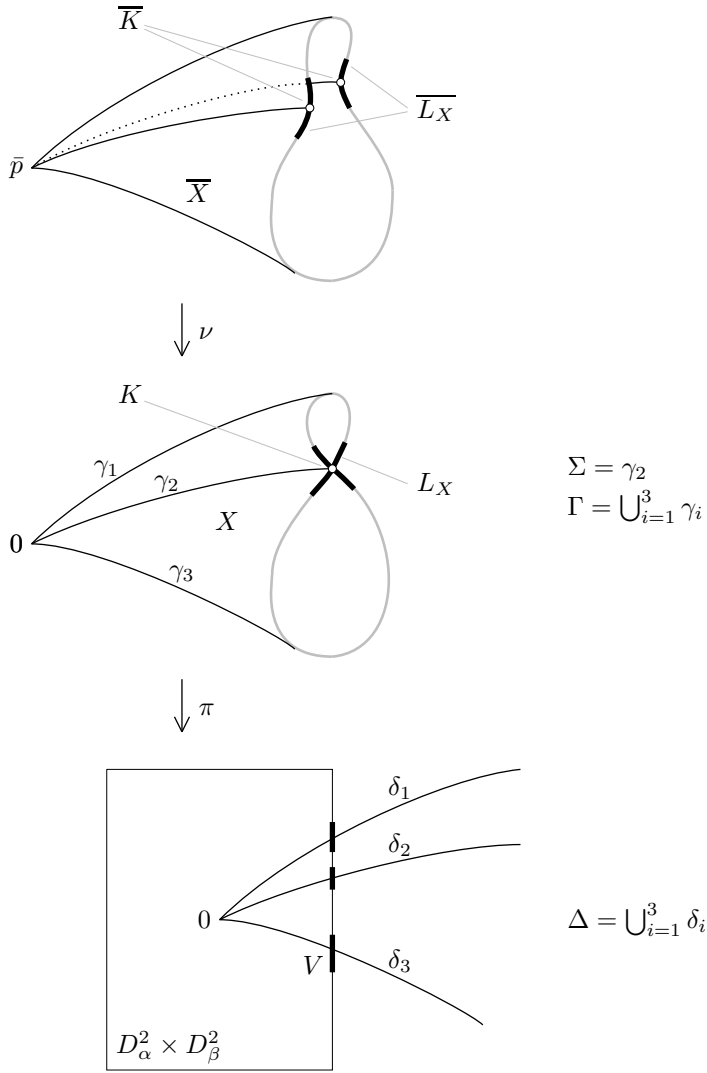

Figure \ref{fig:nor} shows schematically $\bar{\Gamma} =\nu ^{-1} (\Gamma) \subset \bar {X}$  and $\Delta$ when $\Sigma $ is irreducible and $\Gamma\setminus {\Sigma }$ has two  irreducible components.

\begin{lemma} A d-curling is  a  singular pinched solid torus associated to a d-cycle, i.e.  if $c$ is a d-cycle, then  $\cal {C}_d$ is homeomorphic to $T(d(D),c)$.

\end{lemma}

{\it Proof:}  We use the notation of Example 3.4.2.  The model of d-curling obtained  in this example is   the tubular neighbourhood $N(K_{\sigma})$ of the singular knot of the link $L_X$ of $X= \{(x,y,z)\in \C ^3 \ where \  z^d-xy^d=0\}$. As we work up to homeomorphism, it is sufficient to prove that $N(K_{\sigma})$  is a  singular pinched solid torus associated to a d-cycle.
We can  saturate the solid torus $T$ by the oriented circles $l_b= S \times \{b\}, \ b\in  D  $. The circles $\nu (l_b), b\in  D $ also  saturate $N(K_{\sigma})$ with  oriented circles.  The fiber $\pi_x ^{-1}(a)= (C,(a,0,0))$ is a singular fiber of the fibration  $\pi_x : \nu (T) \to S_{\alpha} $. The equation of the curve germ $C$ at $(a,0,0)$ is $ \{ z^d-ay^d=0\}$, this is a plane curve germ with $d$ branches.  So, $C$ is homeomorphic to a $d$-pinched disc. Moreover, the first return along the circles $\nu (l_b)$ is a monodromy  $h$ of $\pi_x$ which satisfies the conditions given in  Definition 3.4.1 to obtain a  singular pinched solid torus associated to a d-cycle.

 Indeed,  $(\pi_x \circ \nu): T \to S_{\alpha} $ is a trivial fibration  with  fiber $ \nu ^{-1} ( C ) = \{(\{u_i\} \times D_{\beta }), u_i^d=a\} $ which is the   disjoint union of   $d$ ordered  meridian discs  of $T$. The first return $h_T$ along the  oriented circles $l_b$  is a cyclic permutation of the ordered $d$  meridian discs and $(h_T)^d$  is the identity morphism. Moreover $\nu $ restricted to $T\setminus (S\times \{0\})$ is a homeomorphism.  As $h_T$ is a lifting of $h$,  by $\nu$,  the monodromy $h$ determines $N(K_{\sigma})$ as a   singular pinched solid torus associated to a d-cycle.

{\it End of proof.}

\begin{proposition} Let $\sigma $ be a  branch of the singular locus of $(X,0)$ which has  a reducible hyperplane section. Let  $(C,s)$ be the hyperplane section germ at $s\in \sigma \cap \{x=a\} $.
Let  $\bar {\sigma}_j,1\leq j \leq n,$ be the $n$ irreducible components  of $\nu _L^{-1}(\sigma)$ and let $d_j$ be the order of  of $\nu_L$ restricted to $\bar {\sigma}_j $.  Let $c_j$ be a $d_j$-cycle and let  $c=c_1 \circ c_2  \circ  \dots  \circ c_n$ be the  permutation of $k(\sigma)=d_1+\dots +d_j+\dots +d_n$  
elements which is the composition of the $n$ disjoint cycles $c_j$. A tubular neighbourhood $N(K_{\sigma})$ of $K_{\sigma }$ is  a singular pinched  solid torus  associated to the permutation $c$. Moreover, the restriction of $\nu $ on $\coprod_{1\leq j  \leq n} N(K_{\bar {\sigma}_j } )$ is the composition of two quotients: the quotients which define the $d_j$-curlings followed  by the quotient $f_{\sigma }$ which identifies their  cores.
 \end{proposition}

{Proof:} Let  $N(K_{\bar {\sigma}_j }), 1\leq j \leq n$  be the $n$ connected components of   $\nu ^{-1}(N(K_{ \sigma }))$. So, $N(K_ {\sigma}) \setminus K_{\sigma}$ has also $n$ connected components and   $(N(K_ {\sigma} ))_j =\nu (N(K_{\bar {\sigma}_j }))$ is the closure of one of them.  Every  $N(K_{\bar {\sigma}_j })$ is a  solid torus  and the restriction of $\nu$ on its  core  $K_{\bar {\sigma}_j }$,  has order $d_j$.  The intersection   $(\nu ^{-1} ( C ) )\cap  N(K_{\bar {\sigma}_j })$ is a disjoint union of $d_j$ ordered and oriented meridian discs of $N(K_{\bar {\sigma}_j })$. We can choose a homeomorphism $g_j: (S\times D) \to N(K_{\bar {\sigma}_j })$ such that $(\nu \circ g_j)^{-1} (  C)=\{u\}\times D,\ u^{d_j}=1 .$ 

The model of a $d_j$-curling gives the quotient  $q_j: (S\times D) \to \cal {C}_{d_j}$.  As in example 3.4.2, there exists  a unique homeomorphism $f_j:  \cal {C}_{d_j} \to (N(K_ {\sigma} ))_j $ such that $f_j\circ q_j=\nu \circ g_j$. So, $(N(K_ {\sigma} ))_j$  is a $d_j$-curling. In particular, if  $\nu _j$ is the restriction of $\nu $ on $N(K_{\bar {\sigma}_j })$, then $\nu_j =f_j \circ q_j \circ (g_j)^{-1}$. Up to homeomorphism $\nu_j$ is equivalent to the quotient which defines the $d_j$-curling.

But  for all $j, 1\leq j \leq n$, we have $\nu  (K_{\bar {\sigma}_j })=(K_{\sigma})$. Up to homeomorphism, $N(K_{\sigma })$ is obtained as the quotient of the disjoint union of the  $d_j$-curlings by the identification of their cores.  The disjoint union of the $f_j$ induces a homeomorphism  $f_{\sigma }$ from $$ N =(\coprod_{1\leq j  \leq n}   \cal {C}_{d_j} )/ q_j(u,0)\sim q_i(u,0) \Leftrightarrow \nu (g_j(u,0))=\nu (g_i (u,0))$$
onto $N(K_ {\sigma} )$. Up to homeomorphism, the restriction of $\nu $ on $\coprod_{1\leq j  \leq n} N(K_{\bar {\sigma}_j } )$ is the composition of two quotients: the quotients which define the $d_j$-curlings followed  by the quotient $f_{\sigma }$ which identifies their cores.
It is sufficient to prove that $N=T(k(\sigma)(D),c)$ where $c$ is the composition of $n$ disjoint cycles $c_j$ of order $d_j$. 
By Lemma 3.4.3,   $\cal {C}_{d_j}=T(d_j(D),c_j)$ and it is obvious that the identifications correspond to the disjoint union of the cycles.

 {\it End of proof.}

\newpage

 \section{  Hirzebruch-Jung's resolution  of $(X,0)$}

  In this Section $(X,0)$ is a normal surface germ.
 
 Let  $\pi  : (X,0)\longrightarrow (\C ^2,0)$  be a finite analytic morphism which is defined  on  $(X,0).$
 For example    $\pi$  can be  the    restriction on $(X,0)$ of a linear projection,  as chosen  in the beginning of Section 3.1.  But the construction can be performed with any finite morphism $\pi$. We denote by    $\Gamma $   the singular locus of $\pi$  and by  $\Delta =\pi (\Gamma )$ its discriminant.

Let   $r: (Z,E_Z)\to (\C^2, 0)$ be the minimal embedded resolution of $ \Delta  $,     $E_Z=r^{-1} (0)$  is  the exceptional divisor of $r$,  and  $E_Z^+=r^{-1} (\Delta )$ is  the total transform of $\Delta.$   The irreducible components of $E_Z$ are smooth because the resolution $r$ is obtained by a sequence of blowing up of points in a smooth surface. Let us denote by $E^0_Z$ the set of the smooth points of $E_Z^+$. So, $E_Z^+ \setminus E^0_Z$ is the set of the double points of  $E_Z^+$.

 Here, we  construct in details  the   Hirzebruch-Jung resolution    $\rho  : (Y,E_Y) \to  (X, 0)$  associated to $\pi.$ This  will prove the existence of a good resolution of $(X,0) $.  As the link $L_X$ is diffeomorphic to the boundary of $Y$, this   will describe $L_X$ as  the boundary of  a plumbing.
 In particular,  we will explain how  to obtain  the dual graph $G(Y)$  of $E_Y$ when we have the dual graph  $G( Z)$ associated to  $E_Z.$    Knowing the Puiseux expansions of all the branches of $\Delta ,$ there exists an algorithm  to compute the dual graph $G_w(Z)$ weighted by the self intersection of the irreducible component of $E(Z)$
 (For example see \cite{b-k} and  Chap. 6 and 7 in \cite{m-w}).  Except in particular cases,  the determination    of  the self intersection of the irreducible components of $E_Y$ is rather delicate.

 \subsection{ First step: Normalization}
 
  We begin with the minimal resolution $r$ of $ \Delta  .$  The pull-back of    $\pi $ by $r $ is a finite morphism  $\pi_r: ( Z ', E_{Z'})\to (Z,E_Z) $  which induces an isomorphism from $E_{Z'}$ to $E_Z$.
 We denote      $r_\pi:  ( Z ', E_{Z'})\to (X,0)$, the pull-back of $r$ by $\pi .$ Figure 3  represents  the  obtained commutative diagram.

 \begin{figure}[h]
$$
\unitlength=0.80mm
\begin{picture}(60.00,60.00)
\thicklines

\put(5,5){\vector(1,0){40}}
\put(5,55){\vector(1,0){40}}
\put(0,50){\vector(0,-1){40}}
\put(50,50){\vector(0,-1){40}}
\put(-12,4){$(Z, E_Z)$}
\put(48,4){$(\C^2,0)$}
\put(48,54){$(X,0)$}
\put(-15,54 ){$(Z',E_{Z'})$}
\put(55,30){$\pi$}
\put(-7,30){$\pi_r$}
\put(25,8){$r$}
\put(25,58){$r_\pi$}
\end{picture}
$$

\caption{ The diagram of the pull-back of the resolution $r$ by $\pi.$
} 
\end{figure}

 In general $ Z'$ is not normal. 
  Let $n  : (\bar Z, E_{\bar Z})\to (Z',E_{Z'})$ be the  normalization of $Z'$.

  \begin{remark}
  \begin{enumerate}
    \item By construction, the discriminant locus of  $\pi_r \circ n$   is included in $ E_Z^+=r^{-1}(\Delta )$ which is the total transform of  $\Delta $ in $Z$. As, $X$ is normal at $0$,  $(X\setminus \{0\})$  has no singular point.
    
  \item As the restriction of $r$ to  $Z\setminus E_Z$ is an isomorphism, the restriction of $r_{\pi}$ to $Z'\setminus E_{Z'}$ is also an isomorphism. Let us denote $\Gamma '$ (resp. $\bar{\Gamma} $) the closure of $(r_{\pi})^{-1}(\Gamma \setminus \{0\})$ in $E_{Z'}$ (resp.  the  closure of $(r_{\pi} \circ n)^{-1}(\Gamma \setminus \{0\})$ in $E_{\bar{Z}}$). The restriction of $r_{\pi}$  on $\Gamma'$ (resp. $(r_{\pi} \circ n)$ on $\bar{\Gamma }$) is  an isomorphism onto  $\Gamma .$
  
  \item The singular locus of $Z'$ is included in $E_{Z'}.$ The normalization $n$ restricted  to $\bar{Z} \setminus  E_{\bar Z}$  is an isomorphism. 
    \end{enumerate}
  \end{remark}

  {\bf {Notations}} We  use the following notations:
   
  $E_{Z'}^+= E_{Z'}\cup \Gamma',$  and   $E_{Z'}^0$ is  the  set of the points of $E_{Z'}$  which belong to a unique irreducible component of   $E_{Z'}^+.$ 
  Similarly:  $E_{\bar {Z}}^+= E_{\bar{Z}}\cup  \bar{\Gamma},$ and   $E_{\bar{Z}}^0$ is  the  set of the points of $E_{\bar{Z}}$   which  belong to a unique irreducible component of  $E_{\bar{Z}}^+.$

  \begin{proposition}
  The only possible singular points of $\bar{Z}$ are points  of  $E_{\bar{Z}}$ which belong to several irreducible components of $E_{\bar{Z}}^+.$    The restriction of the map $(\pi_r \circ n) $
  to   $E_{\bar Z}$  induces a finite morphism  from $E_{\bar Z}$ to $E_Z$ which is a regular covering  from  $( \pi_r \circ n)^{-1}(E_ Z^0)$ to  $(E_ Z^0)$.
  
  \end{proposition}

   {\it Proof:}  As, $X$ is normal at $0$,  $(X\setminus \{0\})$  has no singular point.  The pull-back construction implies that:
   \begin{enumerate}
   \item  The morphism $ \pi_r $  is finite and its generic order is equal to the generic order of $\pi.$ Indeed,  $ \pi_r $  restricted   to $E_{Z'}$ is an isomorphism. Moreover, the  restriction  of   $\pi_r $ to $(Z' \setminus E_{Z'})$ is isomorphic, as a ramified covering,  to the restriction of $\pi $ to $(X\setminus \{0\}).$ So,  the restriction morphism $ (\pi_r)_{\vert } : (Z' \setminus E_{Z'}) \to (Z \setminus E_{Z})$  is   a  finite ramified covering with ramification  locus $\Gamma '.$

\item As the restriction of $r$ to  $(Z\setminus E_Z)$ is an isomorphism,  then   the restriction of $r_{\pi}$ to  $(Z'\setminus E_{Z'})$ is also an isomorphism. 
 So, the restriction of $(r_{\pi}\circ n)$ to  $( \bar{Z}\setminus E_{\bar{Z}})$ is an analytic isomorphism onto the non singular analytic set  $(X\setminus \{0\}).$ It implies that  $( \bar{Z}\setminus E_{\bar{Z}})$ is smooth.
      
       \end{enumerate}

     If  $\bar {P} \in E_{\bar{Z}}^0,$  then   $P=( \pi_r \circ n) ( \bar {P})$ is a smooth point of an irreducible component $E_i$ of  $E_Z.$  The normal fiber bundle to  $E_i$ in $Z$ can be locally  trivialized  at $P$. We can choose a closed small neighborhood $N$ of $P$  in $Z$ such that $ N= D  \times \Delta   $ where  $D$ and $\Delta $ are two discs , $N \cap E_Z =(D \times 0)$ and for all $z\in  D$, $z \times \Delta $ are fibers of the bundle in discs associated to the normal bundle of  basis $E_i$. We choose  $\bar{N} = ( \pi_r \circ n) ^{-1} ( {N})$  as 
 closed  neighborhood  of $\bar{P}$ in $\bar{Z}$.  But $\bar{Z}$ is normal and the local discriminant of the restriction   $(\pi_r \circ n)_{\vert} : (\bar{N},\bar{P}) \to (N,P)$ is included in $D \times 0$ which is a smooth germ of curve. In that case, the link of $(\bar{N} , \bar{P})$ is $S^3$ (in  Lemma 5.2.1,  
 we give   a topological proof of this classical result).  As $\bar{Z}$ is normal,  by Mumford's Theorem \cite{mu},   $\bar{P}$ is a smooth point of $\bar{Z}$. This ends  the proof of the first statement of the proposition.
 
 Now,  we know that  the morphism  $(\pi_r \circ n)_{\vert} : (\bar{N},\bar{P}) \to (N,P)$ is a finite morphism between two smooth germs of surfaces  with  non singular discriminant locus. Let $d$
be its generic order.   
 By Lemma 5.2.1,  such a morphism is locally isomorphic (as an  analytic morphism) to  the morphism defined on $(\C^2,0)$ by   $(x,y)  \mapsto  (x,y^d).$ So,   $\bar{D} = ( \pi_r \circ n) ^{-1} ( {D \times 0})$ is a smooth disc in $E_{\bar{Z}}^0$ and 
 the restriction of such  a morphism to $\{ (x,0), x\in \bar{D}\}$ is a local isomorphism. 
 
 By definition of $E^0_Z$, $P\in (E_i\cap E^0_Z)$ is a smooth point in the total transform of $\Delta.$ If we take a smooth germ $(\gamma, P)$ transverse to $E_i$ at $P$, then  $(r(\gamma ), 0)$ is not a branch of $\Delta .$ The restriction of $\pi$ to $ \pi ^{-1} (r(\gamma ) \setminus 0)$ is a regular covering. Let $k$ be the number of irreducible components of $ \pi ^{-1}( r(\gamma ))$. The number  $k$ is constant for all $P\in E_i \cap E^0_Z .$ Let $P'$ be the only point of $ (\pi _r )^{-1} (P) $. Remark 3.3.5 which uses Lemma 3.3.4,  shows that the $k$ irreducible components of the  germ of curve $((\pi_r )^{-1}(\gamma),P') $  are separated by $n.$ So,  the restriction of the map $(\pi_r \circ n)$ to  $(( \pi_r \circ n)^{-1}( E_i \cap E_ Z^0))$ is a regular covering of order $k.$

    {\it End of proof.}

  \begin{definition} 
   A germ $(W,0)$ of  complex  surface is {\bf  quasi-ordinary} if there exists a finite morphism 
  $\phi :(W,p) \to (\C^2,0)$ which has  a  normal-crossing discriminant.
 A {\bf  Hirzebruch-Jung singularity }  is a quasi-ordinary singularity of normal  surface germ.

\end{definition}

    \begin{lemma}  Let  $\bar {P}$ be  a point of $E_{\bar{Z}}$ which belongs to several irreducible components of $E_{\bar{Z}}^+.$ Then  $\bar {P}$ belongs to two  irreducible components of $E_{\bar{Z}}^+.$ Moreover,   either  $\bar {P}$ is a smooth point of $\bar{Z}$ and  $E_{\bar{Z}}^+$ is a normal crossing divisor around $\bar {P},$ either  $\bar P$ is a Hirzebruch-Jung singularity of $\bar Z$.  
 \end{lemma}

 {\it Proof:}  If  $\bar {P}$ be  a point of $E_{\bar{Z}}$ which belongs to several irreducible components of $E_{\bar{Z}}^+$    then,  $P=( \pi_r \circ n) ( \bar {P})$  is a double point of $E_Z^+$.  Moreover $Z$ is smooth and $E_Z^+$ is a normal crossing divisor. We can choose  a closed neighbourhood  $N$ of $P$  isomorphic to a product of disc $(D_1\times D_2),$ and we take  $\bar{N} = ( \pi_r \circ n) ^{-1} ( {N})$. For a sufficiently small $N,$ the restriction of $(\pi_r \circ n)$ to  the pair $(\bar {N}, \bar {N} \cap E_{\bar{Z}}^+)$  is a finite ramified morphism  over the pair   $(\bar {N}, \bar {N} \cap E_{\bar{Z}}^+)$  and the ramification locus is included in the normal crossing divisor   $(N \cap E_Z^+ ) .$  The pair $(\bar {N}, \bar {P})$ is normal and the link of the pair $( {N},  {N} \cap E_Z^+)$ is the Hopf link in $S^3.$
  Then the link of $\bar {N}$ is a lens space, and  the link of   $( \pi_r \circ n) ^{-1}({N} \cap E_Z^+) $ has two components (Lemma 5.1.4  gives  a topological proof of this classical result). So,  $E_{\bar{Z}}^+$ has two irreducible components at $\bar {P}$. We have two possibilities:
  \begin{enumerate}
  \item $\bar{P}$ is a smooth point  in $\bar{Z}$. Then  the link of the pair  $(\bar {N}, \bar {N} \cap E_{\bar{Z}}^+)$  is the Hopf link in $S^3$ and $E_{\bar{Z}}^+$  is a normal crossing divisor at $\bar{P}$. 
  
  \item   $\bar{P}$ is an isolated singular point of $\bar{Z}$. Then, the link of $\bar{N}$ is a lens space which is not $S^3$.  The point $\bar{P}$ is a Hirzebruch-Jung singularity of $\bar{Z}$ equipped  with  the finite morphism

 $$ (\pi_r \circ n )_| :  (\bar{N}, \bar{N}\cap E_{\bar{Z}}^+)  \to   (N, N\cap E_{Z}^+)  $$ which has the normal crossing divisor  $N\cap E_{Z}^+$ as discriminant. 
  \end{enumerate}

  {\it End of proof.}

  The example given in  Section 6 illustrates the following Corollary.
  
  \begin{corollary}
   Let  $G(\bar{Z})$ be the dual graph   of $E_{\bar{Z}}$.  Proposition 4.1.2  and Lemma 4.1.4  imply that $(\pi_r \circ n)$ induces a finite ramified covering of graphs from   $G(\bar{Z})$ onto $G(Z)$.
 
 \end{corollary}

\newpage

 \subsection{ Second  step:  Resolution of the Hirzebruch-Jung singularities}

\bigskip
 If   $\bar P$ is a singular point of $\bar Z$, then $P=(\pi_r \circ n) (\bar P)$ is a double point of $E_Z^+$ . In particular, there are finitely many  isolated  singular points in $\bar Z$.
The singularities of $\bar Z$ are Hirzebruch-Jung singularities. More precisely, 
let $\bar P_i, 1\leq i\leq n$, be the finite set of the singular points of $\bar Z$ and let  $\bar{U_i}$  be a sufficiently small neighborhood of  $\bar P_i$ in $\bar Z$. 
We have the following result (see \cite{hi} for a proof, see also \cite{la},   \cite{po}  and  \cite{l-w}) and, to be self-contained, we give a proof  in  Section 5.3 (Proposition 5.3.1):

{\bf Theorem}
The exceptional divisor of the minimal resolution of $(\bar {U_i}, \bar { P_i)} $ is a normal crossings divisor with  smooth rational irreducible components and its    dual graph is a bamboo (it means  is homeomorphic to a segment).

\

\bigskip
\bigskip
\bigskip

 \begin{figure}[h]
$$
\unitlength=0.80mm
\begin{picture}(120.00,150.00)
\thicklines

\put(70,5){\vector(1,0){30}}
\put(70,55){\vector(1,0){30}}
\put(60,50){\vector(0,-1){35}}
\put(110,50){\vector(0,-1){35}}
\put(50,4){$(Z,E_Z)$}
\put(103,4){$(\C^2,0)$}
\put(103,54){$(X,0)$}
\put(50,54 ){$ (Z', E_{Z'})$}
\put(115,30){$\pi$}
\put(50,30){$\pi_r$}
\put(85,8){$r$}
\put(85,58){$r_\pi$}
\put(27,88){\vector(1,-1){25}}
\put(15,92){$(\bar Z, E_{\bar Z})$}
\put(15,165){$(Y, E_{Y})$}
\put(22,160){\vector(0,-1){60}}
\put(35,165){\vector(3,-1){70}}
\put(15,130){$\bar\rho$}
\put(103,135){$( Y', E_{Y'})$}
\put(110,130){\vector(0,-1){65}}
\put(38,80){$n$}
\put(115,100){$\rho '$}
\put(65,160){$\beta$}
\put(30,160){\vector(3,-4){70}}
\put(65,120){$\rho$}
\end{picture}
$$
\caption{The commutative diagram of  the morphisms involved in the Hizebruch-Jung  resolution  $\rho$ of $\pi$.  By  construction  $\rho =  r_\pi \circ n \circ \bar \rho .$ 
} 
\end{figure}
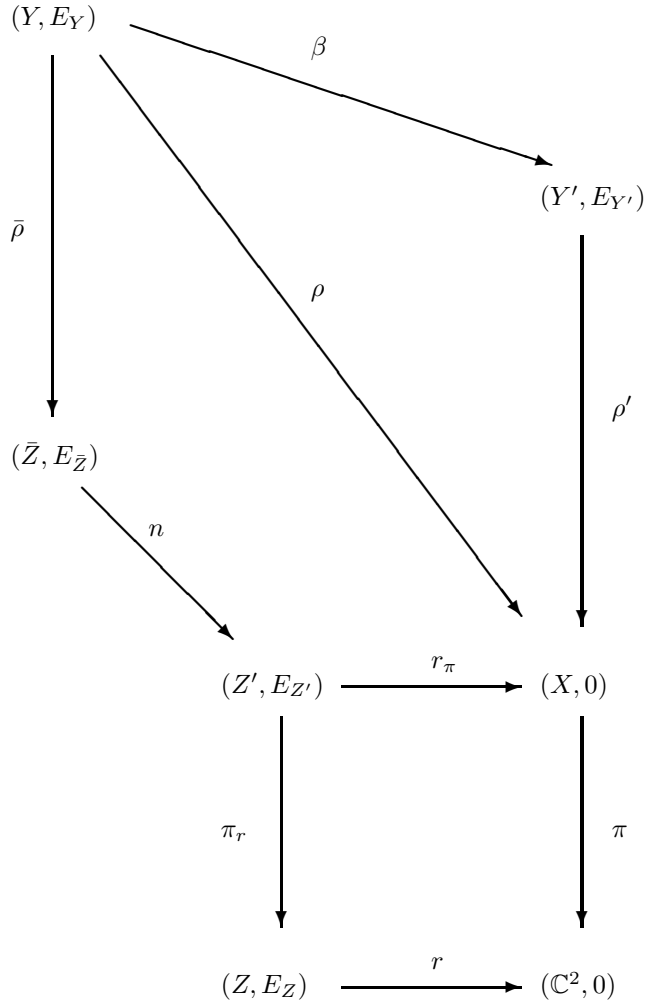

\newpage

Let $\bar \rho _i :(U'_i, E_{U'_i})\to (\bar U_i, \bar {P_i})$ be the minimal  resolution of the singularity $(\bar{U_i}, \bar{(P_i})$. 
From  \cite{l-w} (corollary 1.4.3), see also  \cite{po} (paragraph 4), the spaces $ U'_i$ and the maps $\bar \rho_i$ can be gluing,  for $1\leq i\leq n$,  in  a suitable way  to give a smooth space $Y$ and a map $\bar \rho : (Y, E_Y)  \to (\bar Z, E_{\bar Z})$   satisfying the following property.

\begin{theorem}
Let us denote $\rho=  r_\pi \circ n \circ \bar \rho $. 
Then,  $  \rho :  (Y, E_Y)  \to (X,p)$ is  a good resolution of the singularity $( X,p)$ in which the total  transform  $\rho ^{-1}(\Gamma)=E_Y^+$  of  the singular locus $\Gamma $  of $\pi$ is a normal crossings divisor.
\end{theorem}

{\it Proof:}  The surface $Y$ is smooth because   $\bar{\rho}$  is a resolution of all the singular points of $\bar{Z}.$ As proved in   Proposition 4.1.2 and Lemma 4.1.4,
  the only  possible singular points of the irreducible components of $E_{\bar Z}$ are  the  double points  $\bar {P_i}$ of $E_{\bar Z}^+$ . These points are solved by the resolutions $\bar \rho _i.$ So, the strict transform, by $\bar \rho $, of the irreducible components of $E_{\bar {Z}}$  are smooth. 
  
  The irreducible components of $E_Y$ created during the resolution  $\bar{\rho}$  are smooth rational curves. So, all the irreducible components of $E_Y$ are smooth complex curves.

  By Lemma 4.1.4,  the only possible points of $E_{\bar {Z}}^+$ around  which $E_{\bar {Z}}^+$  is  not smooth or a normal crossing divisor are the Hirzebruch-Jung singularities  $\bar {P_i}, 1\leq i\leq n$. But  as the  $\bar \rho _i, \ 1\leq i\leq n $,   are good resolutions of these singularities, $((\bar \rho _i)^{-1}  (\bar U_i)) \cap (E_Y^+), 1\leq i \leq n,$ are  normal crossing divisors. 
  
   {\it End of proof}.

 As $\rho $ is the composition of three well defined morphisms  which depend only on the choice of the morphism $\pi$ and as we follow the  Hirzebruch-Jung   method, we have the following definition.

\begin{definition} 
 The morphism $ \rho :  (Y, E_Y) \to  (X, 0)$ is the Hirzebruch-Jung   resolution   associated to $\pi$.
\end{definition}

  \begin{corollary}
   The  dual graph $G(Y)$  of $E_Y$ is obtained from the dual graph $G({\bar Z})$ of $E_{\bar Z}$ by replacing the edges,  which represent the Hirzebruch-Jung singular points of $\bar{Z}$,  by a bamboo.
 
 \end{corollary}

 Now we can use the following result (for a proof see   \cite{la}, Theorem 5.9, p.87):

{\bf Theorem} 
 Let  $\rho'  : (Y',E_{Y'})\to (X,0)$ be the minimal resolution of $(X,0)$. There exists $\beta :  (Y, E_Y) \to (\tilde Y',E_{Y'})$ such that $\rho' \circ \beta=\rho$
and the map $\beta$ consists in a composition of  blowing-downs  of  irreducible components,   of the successively obtained exceptional divisors,
of self-intersection $-1$ and   genus $0$,   which are not rupture components.

\newpage

\section{ Appendix:{\bf The topology of a  quasi-ordinary singularity of surface}}

\subsection{Lens spaces}

One can find details on lens spaces and surface singularities in \cite{we}. See also \cite{po1}.

\begin{definition} A lens space $L$ is an oriented compact three dimensional topological manifold  which can be obtained as the union of two solid tori  $T_1\cup T_2$ along their boundaries. The torus $\tau =T_1\cap T_2 $ is  the  Heegaard torus of the  given decomposition $L=T_1\cup T_2.$
   
\end{definition}  
\begin{remark} If $L$ is a lens space, there exists  an embedded torus $\tau $ in $L$ such that $L\setminus \tau $ has two connected components  which are open solid tori $\mathring{T_i} , i=1,2$. Let $T_i, i=1,2,$ be the two compact solid tori  closure of  $\mathring{T_i} $ in $L.$ Of course $\tau =T_1 \cap T_2.$  In \cite{bona}, F. Bonahon show that  a lens space has   a unique, up to isotopy,   Heegaard torus. This  implies that the decomposition $L=T_1 \cup T_2$ is unique up to isotopy, it is  `` the"  Heegaard decomposition of $L$.
\end{remark}

A lens space $L$ with a decomposition of Heegaard torus $\tau $ can be described  as follows. The solid tori $T_i, i=1,2,$ are oriented by the orientation induced by $L$. Let $\tau _i$ be the torus $\tau $ with the orientation induced by $T_i$. By definition a meridian  $m_i$ of $T_i$ is a closed oriented  circle on $\tau_i $ which is  the boundary of a disc $D_i$ embedded in  $T_i$. A meridian of a solid torus is well defined up to isotopy. A parallel $l_i$ of $T_i$  is a closed  oriented curve on $\tau_i$ such that the intersection $m_i\cap l_i=+1$ (we  also write $m_i$ (resp.  $l_i$) for  the homology class of $m_i$ (resp. $l_i$) in the first homology group of $\tau_i $). The homology classes  of two parallels differ by a multiple of the meridian. 

 We choose  on $\tau_2 $, an oriented meridian $m_2$  and a parallel $l_2$ of the solid torus $T_2.$ 
 As in \cite{we}, p. 23, we write a meridian $m_1$ of $T_1$ as  $m_1= nl_2-qm_2$ with $n\in \N$ and $q\in \Z$ where $q$ is well defined modulo $n$.  As $m_1$ is a closed curve on $\tau$, $q$ is prime to $n$. Moreover, the class of  $q$ modulo n  depends on the choice of $l_2.$ So, we can chose $l_2$ such that $0\leq q<n.$ 
 
 \begin{definition}
 
 By a Dehn filling argument,   it is sufficient to know the homology class  $m_1= nl_2-qm_2$ to reconstruct $L.$  By definition {\bf the lens space $L(n,q)$} is the lens space constructed with $m_1= nl_2-qm_2$. We have two very particular cases:
 
 1) $m_1= m_2$, if and only if  $L$ is homeomorphic to $S^1\times S^2$,
 
 2) $m_1= l_2$ if and only if  $L$ is homeomorphic to $S^3.$
 
\end{definition}

\begin{lemma}

 Let  $\phi  :(W,p) \to (\C^2,0)$ be a finite morphism defined  on  
a irreducible  surface germ $(W,p).$ If  the discriminant $\Delta $ of $\phi$ is included in a normal crossing germ of curve,  then the link $L_W$ of $(W,p)$ is a lens space. The link $K_{\Gamma }$  of the singular locus $\Gamma $ of $\phi$,  has at most two connected components. Moreover,  $K_{\Gamma }$    is a sub-link of the  two cores of the two solid tori  of a  Heegaard decomposition of $L_W$ as  a union of two solid tori.
\end{lemma}

 {\it  Proof:}   After performing  a possible analytic  isomorphism of $(\C^2,0)$, $\Delta $ is,  by hypothesis, included in the two axes i.e. $\Delta \subset \{xy=0\}$.

Let $D_{\alpha }\times D_{\beta} \in \C^2$ be a polydisc at the origin  in $\C^2$  where   $0<\alpha <\beta<\epsilon$  are chosen sufficiently small as in Section 3.1. Then, the restriction $\phi _L$ of $\phi$ on the link $L_W$ is a ramified covering on the sphere (with corners)

$${\cal S} = (S_{\alpha }\times D_{\beta}) \cup  (D_{\alpha }\times S_{\beta})$$ with a set of ramification values included in the Hopf link $ K_{xy}= (S_{\alpha } \times \{0\}) \cup  (\{ 0\}  \times S_{\beta})$.

 Let $N(K_{xy})$ be a small  compact tubular neighborhood of $K_{xy}$ in ${\cal S}.$  Then, $N(K_{xy})$ is the union of two disjoint  solid tori  $T_y= (S_{\alpha } \times D_{\beta '}), 0<\beta'<\beta $,   and   $T_x= (D_{\alpha' } \times S_{\beta }), 0<\alpha'<\alpha $.  Then,  $\phi _L^{-1}(T_x)$ (resp.   $\phi _L^{-1}(T_y)$) is an union of $r_x>0$ (resp. $r_y>0$ ) disjoint solid tori because the set of the ramification values of $\phi_L $ is included in the core of $T_x$ (resp. $T_y$).

  Let $V$ be the closure, in ${\cal S}, $ of ${\cal S} \setminus  N(K_{xy}).$ But, $V$ is a thickened torus which does not meet the ramification values of $\phi _L$.  Then,  $\phi _L^{-1}(V)$ is a  union of $r>0$ disjoint thickened tori. But, $L_W$ is connected because $(W,p)$ is  irreducible by hypothesis. The only possibility to obtain  a connected space by gluing  $\phi _L^{-1}(T_x)$,   $\phi _L^{-1}(T_y)$ and   $\phi _L^{-1}(V)$ along their boundaries is  $1=r=r_x=r_y.$

  So,   $\phi _L^{-1}(T_x)$ (resp.   $\phi _L^{-1}(T_y)$) which is in $L_W$ a deformation retract of $T_2=  \phi _L^{-1}( S_{\alpha } \times D_{\beta })$ (resp. $T_1=  \phi_L^{-1}(D_{\alpha } \times S_{\beta })$ ) is  a  unique solid torus. Then $ \tau = \phi _L^{-1}(S_{\alpha} \times S_{\beta})$ is a unique torus.
  We have proved that $L_W$ is the lens space  obtained  as the union of the two solid tori $T_1$ and $T_2$ along their common  boundary $ \tau = \phi _L^{-1}(S_{\alpha} \times S_{\beta})$.  So, $T_1\cup T_2$ is  a  Heegaard decomposition of $L_W$ as  a union of two solid tori.

  By hypothesis $ K_{\Delta } \subset (S_{\alpha } \times \{0\}) \cup  (\{ 0\}  \times S_{\beta})$. Then,  $K_{\Gamma}$ is included in the disjoint union  of  $ \phi _L^{-1}(S_{\alpha } \times \{0\}) $ and $ \phi _L^{-1} (\{ 0\}  \times S_{\beta})$ which are  the cores of  $T_1$ and $T_2$. So, $K_{\Gamma }$ has at most two connected components.

   {\it End of proof.}
   
   \begin{example} Let $n$ and $q$ be two positive integers prime to each other. We suppose that  $0<q<n$. Let $X=\{(x,y,z)\in \C ^3  \ s.t.\  z^n-xy^q=0\}$. The link $L_X$ of $(X,0)$ is the lens space $L(n,n-q)$. 
   \end{example}
   Indeed, let $\phi : (X,0)\to (\C^2,0)$ be the projection  $(x,y,z) \mapsto  (x,y)$ restricted to $X$.
   The discriminant  $\Delta $ of $\phi$ is equal to $ \{xy=0\}$. By Lemma 5.1.4,   $L_X$ is a lens space. We use the notations of 5.1.4,  $L_X= \pi ^{-1}(\cal{S})$ where
   $${\cal S} = (S_{\alpha }\times D_{\beta}) \cup  (D_{\alpha }\times S_{\beta}).$$ 
 In the proof of Lemma 5.1.4, it is shown that  $T_2=  \phi ^{-1}( S_{\alpha } \times D_{\beta })$  and  $T_1=  \phi^{-1}(D_{\alpha } \times S_{\beta })$ )  are two solid tori.
 Let $(a,b)\in   (S_{\alpha }\times S_{\beta})$.  As $n$ and $q$ are prime to each other   $ m_1= \phi ^{-1}(\{a\} \times S_{\beta})$  and  $ m_2 = \phi ^{-1}(S_{\alpha} \times \{b\})$ are connected. So, $m_i,\  i=1,2$,  is a meridian of $T_i$. 
 
 We choose $c\in \C$  such that $c^n=ab^q$. Let $l_2=  \{z=c\} \cap \phi ^{-1}(S_{\alpha} \times S_{\beta})$. On the torus  $ \tau = \phi ^{-1}(S_{\alpha} \times S_{\beta})$, oriented as the boundary of $T_2$, we have  $m_2\cap l_2=+1$ and    $m_1= nl_2-(-q)m_2$. As defined in 5.1.3, we have  $L_X=L(n,-q)=L(n, n-q)$.

\subsection{Finite morphisms with smooth discriminant}

\begin{lemma} 

  Let  $\phi  :(W,p) \to (\C^2,0)$  be a finite morphism,   of generic order $n$,  defined  on  
a normal surface germ $(W,p).$   If the discriminant of $\phi $ is a smooth germ of curve, then     $(X,0)$ is analytically isomorphic to $(\C^2,0)$ and   $\phi $ is analytically  isomorphic to  the map from $(\C^2,0)$  to $(\C^2,0)$ defined by $(x,y)  \mapsto (x,y^n).$

\end{lemma}

{\it  Proof:}  After performing  an  analytic  automorphism  of $(\C^2,0)$, we can  choose coordinates such that $\Delta=\{y=0\} $.

Let $D_{\alpha }\times D_{\beta} \in \C^2$ be a polydisc at the origin  in $\C^2$  where   $0<\alpha <\beta<\epsilon$  are chosen sufficiently small as in Section 3.1. Then, the restriction $\phi _L$ of $\phi $ on the link $L_W$ is a ramified covering on the sphere (with corners) $${\cal S} = (S_{\alpha }\times D_{\beta}) \cup  (D_{\alpha }\times S_{\beta})$$ with a set of ramification values included in the trivial  link $ K_{y}= (S_{\alpha } \times \{0\})$.

 Here,  we satisfy the hypothesis  of Lemma 5.4. So, $T_2=  \phi _L^{-1}( S_{\alpha } \times D_{\beta })$ and $T_1=  \phi_L^{-1}(D_{\alpha } \times S_{\beta })$  are two  solid tori with common boundary  $ \tau = \phi _L^{-1}(S_{\alpha} \times S_{\beta})$.  We take $a\in S_{\alpha }$ and $b\in S_{\beta }.$  

Let us consider ${\cal D}_a=  \phi _L^{-1}( \{a\} \times D_{\beta }) \subset T_2$ and ${\cal D}_b=  \phi_L^{-1}(D_{\alpha } \times \{b\}) \subset T_1.$ 

Here the singular locus of $\phi _L$ is the core of $T_2$ and does not meet $T_1$. 

The restriction of $\phi _L$ on $ \phi _L^{-1}( a \times b)$is a regular covering over a disc. Then ${\cal D}_b$ is a disjoint union of $n$ discs where $n$ is the general order of $\phi_L.$  Let $m_1$ be the oriented boundary of one of the $n$   discs which are the  connected components of  ${\cal D}_b$. By definition $m_1$ is a meridian of $T_1.$

The restriction of $\phi _L$ on ${\cal D}_a$ is a  covering over a disc and   $(a\times 0)$ is  the only  ramification value . Then ${\cal D}_a$ is a disjoint union of $d$ discs where $d<n$. On $\tau $,  the intersection between  the circles boundaries of ${\cal D}_a$ and ${\cal D}_b$ is  equal to $n$ because it is given by the (positively counted) $n $ points  of $ \phi _L^{-1}( a \times b)$.  The restriction of $\phi_L$ on $T_1$ is a Galois covering of order $n$ which permute cyclically the connected components of ${\cal D}_a$. So, on the torus $\tau =b(T_1)$,  any of the $d$ circles boundaries  of the  connected components of 
${\cal D}_a$ cuts any of the $n$  circles  boundaries of the connected components  of ${\cal D}_b$. So computed, the intersection $b({\cal D}_a)\cap b({\cal D}_b)$ is equal to $nd$. But, $nd=n$ because  this intersection is given by the $n$ points of  $ \phi _L^{-1}( a \times b)$.

So,  $d=1$ and ${\cal D}_a$ has a unique connected component. The boundary of ${\cal D}_a$ is a meridian $m_2$ of $T_2$. As $m_1$ is the boundary of one of the $n$ connected components of ${\cal D}_b$, $m_1 \cap m_2 =+1$ and $m_1$ can be a parallel  $l_2$ of $T_2$. This  is the case 2) in Definition 5.1.3, so  the link $L_W$ of $(W,p)$ is  the 3-sphere $S^3$. 
 As $(W,p)$ is normal, by Mumford \cite{mu}, $(W,p)$ is a smooth surface germ i.e  $( W,p)$ is  analytically  isomorphic to $(\C^2,0)$. The first part of Lemma 5.5 has been proved.

(*) Moreover $\phi_L^{-1} (S_{\alpha } \times \{0\}) \cup  (\{ 0\}  \times S_{\beta})$ is the union of the cores of $T_1$ and $T_2$. Then, {\bf  $(S_{\alpha } \times \{0\}) \cup  (\{ 0\}  \times S_{\beta})$ is a Hopf link in the 3-sphere $L_W$}.

From now on, $\phi  :(\C^2,0) \to (\C^2,0)$ is a finite morphism and its discriminant locus is $\{y=0\}$.
Let us write  $\phi=(\phi_1, \phi_2)$. The zero locus of the function germ $$(\phi_1.\phi_2)  :(\C^2,0) \to (\C^,0)$$ is the link  describe just above (see (*)), i.e.  it is   a Hopf link.  The  function $(\phi_1 .\phi_2) $ reduced is analytically  isomorphic to $(x,y)  \mapsto (xy).$ But $\phi_1 $ is reduced because its  Milnor fiber  is diffeomorphic to  ${\cal D}_a=  \phi _L^{-1}( \{a\} \times D_{\beta }) \subset T_2$  which is a disc. So, $\phi_1$ is isomorphic to  $x$. 

The Milnor fiber of $\phi_2 $ is diffeomorphic to the disjoint union of the  $n$  discs  ${\cal D}_b=  \phi_L^{-1}(D_{\alpha } \times \{b\}) \subset T_1.$ When the Milnor fiber of a function germ $f: (\C^2,0) \to (\C,0)$ has $n$ connected component,  $n$ is the $g.c.d.$ of the multiplicity of the irreducible factors of $f$. Here $\phi_2=g^n$ where $g$ is an irreducible function germ. We already  have seen that $\phi_2$ reduced is isomorphic to $y$. This achieves  the proof that $\phi_2$ is isomorphic to $y^n$ and  $\phi=(\phi_1, \phi_2)$ is isomorphic to $(x,y^n)$.

{\it End of proof}.

\newpage

\subsection{The Hirzebruch-Jung singularities}

 \begin{proposition}  Let   $(W,p)$ be a normal surface germ such that there exists  a finite morphism $\phi  :(W,p) \to (\C^2,0)$ which has a normal-crossing discriminant $(\Delta ,0)$. Then, $(W,p)$ has a minimal good resolution $\rho : (\tilde W, E_{\tilde W}) \to (W,p)$ such that :
  \item I)  the exceptional divisor $E_{\tilde W}$ of $\rho$ has smooth rational  irreducible components and its  dual graph  is a bamboo. We orient the bamboo from the vertex (1) to the vertex $(k)$. The vertices are indexed by this orientation, 

\item II) the strict  transform   of $\Phi ^{-1} ( \Delta )$  has two smooth  irreducible components which meet transversaly   $E_{\tilde W}$, one of them  at a smooth point of $E_1$ and  the other  component at a smooth point of $E_k$.
  \end{proposition}
  
  {\it  Proof:}  After performing  an  analytic  isomorphism of $(\C^2,0)$, we can  choose coordinates such that $\Delta=\{xy=0\} $. We have to prove that there exists a minimal resolution $\rho $ of $(W,p)$  such that   the shape of the dual graph of the total transform of $\Delta $ in  $\tilde W$  looks like   the graph drawn in   Figure 5 where  all vertices represent smooth rational curves.

  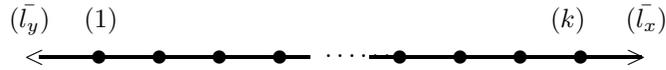
\begin{figure}[h]
$$
\unitlength=0.40mm
\begin{picture}(40.00,20.00)
\thicklines
\put(-110,10){$(\bar{l_y})$}
\put(95,10){$(\bar {l_x})$}
\put(-85,10){$(1)$}
\put(70,10){$(k)$}
\put(-100,0){\line(1,0){90}}
\put(10,0){\line(1,0){90}}
\put(-105,-2.3){$<$}
\put(95,-2.3){$>$}
\put(-80,0){\circle*{4}}
\put(-60,0){\circle*{4}}
\put(-40,0){\circle*{4}}
\put(-20,0){\circle*{4}}
\put(-5,0){\ldots\ldots }
\put(20,0){\circle*{4}}
\put(80,0){\circle*{4}}
\put(60,0){\circle*{4}}
\put(40,0){\circle*{4}}

\end{picture}
$$

\caption{ The shape of the dual graph of  $G(\tilde W)$  to which  we add an arrow  to  the vertex $(1)$   to represent the strict transform  of $\{ x=0\}$ and another   arrow  to the vertex $(k)$ to represent the strict transform of $\{ y=0\}$. } 

\end{figure}

By lemma 5.1.4, the link  $L_W $ of $(W,p)$ is a lens space. If $L_W$ is homeomorphic to $S^3$, $(W,p)$ is smooth by Mumford \cite{mu},  and there is nothing to prove. Otherwise,  let $n$ and $q$ be  the two positive integers, prime to each other, with $0<q<n$,  such that $L_W$ is the  lens  space $L(n,n-q)$. By Brieskorn \cite{br1} (see also Section 5 in  \cite{we}),  the normal quasi-ordinary complex surface germs  are taut. It  means that any normal quasi-ordinary complex surface germ $(W',p')$ which has a link orientation preserving  homeomorphic to $L(n,n-q)$ is  analytically  isomorphic  to $(W,p)$. In particular, $(W,p)$ and $(W', p')$ have isomorphic minimal good resolutions. Now, it is sufficient to describe the good minimal resolution of a  given normal quasi-ordinary surface  germ which has  a link homeomorphic to  $L(n,n-q)$. As explained  below, we can use $({\bar X}, {\bar p})$  where  $\nu : ({\bar X}, {\bar p})  \to  (X,0)$  is the normalization of  $X=\{(x,y,z)\in \C ^3  \ s.t.\  z^n-xy^q=0\}$.

   \begin{lemma} Let $n$ and $q$ be two positive integers prime to each other. We suppose that  $0<q<n$. Let $X=\{(x,y,z)\in \C ^3  \ s.t.\  z^n-xy^q=0\}$.  There exists a good resolution $\rho_Y : (Y, E_Y) \to (X,0)$ of $(X,0)$ such that the dual graph $G(Y)$ of $E_Y$  is a bamboo and the dual graph of the total transform of $\{xy=0\}$ has the shape of the graph given in Figure 5.
   \end{lemma}
     Lemma 5.3.2 implies Proposition 5.3.1. Indeed:
     
   1)  In  Example 5.1.5,  we show that the link $L_X$  of $(X,0)$ is the lens space $L(n, n-q)$.  Let    $\nu : ({\bar X}, {\bar p})  \to  (X,0)$ be the normalization of $(X,0)$. The singular locus of $(X,0)$ is the line $\Sigma =\{ (x,0,0), x\in \C \}$.  For $a\in \C$, the hyperplane section of $X$ at $(a,0,0)$ is the  plane curve germ $\{z^n-ay^q=0\}$. As $n$ and $q$ are prime to each other $\{z^n-ay^q=0\} $ is irreducible. Then,  by  Proposition 3.3.6, $\nu $ is a homeomorphism. So,  the link $L_{\bar X}$  of $({\bar X},0)$ is the lens space $L(n, n-q)$.
   
   2) Let $\rho_Y : ( Y, E_Y) \to (X,p)$  be a good resolution of $(X,0)$  given as in  Lemma 5.3.2, in particular  the dual graph  $G(Y)$  of $E_Y$ is a bamboo. As  any  good resolution factorizes   by the normalization  $\nu : (\bar X,\bar p) \to (X,0 )$ (see \cite{la} Thm. 3.14), there exists  a unique  morphism $\bar{\rho_Y }: ( Y, E_Y) \to (\bar{X},\bar {p})$  which is a good resolution of $(\bar{X},\bar {p})$ .  Let $ \rho' : ({ Y'}, E_{ Y'})\to (\bar {X},\bar {p})$ be the   minimal good  resolution of  $(\bar {X},\bar {p})$. 
   Then,   (for example see \cite{la} Thm 5.9 or \cite{b-p-v} p. 86), 
 there exists  a morphism $\beta : (Y,E_Y)\to  (Y',E_{Y'})$  which is a  sequence of  blowing-downs of  irreducible components of genus zero and self-intersection $-1$.  By Lemma 5.3.2,  the dual graph   $G(Y)$  is a bamboo and the dual graph of the total transform of $\{xy=0\}$ has the shape of the graph given in Figure 5. So, the morphism of graph $\beta *: G(Y) \to G(Y')$ induced by $\beta $, is  only   a  contraction  of $G(Y)$ in a shorter bamboo.
 
 {\it  Proof of Lemma 5.3.2}:
 
 In $X$, we consider  the lines $l_x=\{(x,0,0), x\in \C\}$ and $l_y=\{(0,y,0),y\in \C\}$ and  the singular locus of $(X,0)$ is equal to $l_x$. We prove Lemma 5.3.2 by a finite induction on $q\geq 1$.
 
 1) If $q=1$, $X=\{(x,y,z)\in \C ^3  \ s.t.\  z^n-xy=0\}$ is the well-known normal singularity $A_{n-1}$. The minimal resolution is a bamboo of $(n-1)$ irreducible components of genus zero. Indeed,  to construct $\rho_Y : (Y, E_Y) \to (X,0)$,  it is sufficient to perform  a sequence  of  blowing-ups  of   points ( we blow up $n/2$ points when $n$ is even and  $(n-1)/2$ points  when $n$ is odd) . We begin to  blow up the origin, this  separates the strict transform of the  lines $l_x$ and $l_y$. The exceptional divisor,    in the strict transform of $(X,0)$ by  the blowing-up of    the origin in $\C^3$,   has two irreducible rational components when  $n>2$ and only one  irreducible rational component  when $n=2$. If $n>2$, we continue by  the  blowing-up  of  the intersection point of the two irreducible components of the exceptional divisor.
 
 2)  Let us state how we proceed if $1<q<n$:
 
 I)  As $n$ and $q$ are prime to each other, the rest $r$ of the division $n=mq+r$ is prime to $q$ and $1<r<q$. Let  $R: Z \to \C^3$ be a sequence of $m$ blowing-ups of the line $l_x$ in $\C^3$  and of  its strict transforms in a smooth  3-dimensional complex space.  Let $Y_1$ be the strict transform of $X$ by $R$.  Let  $\rho : (Y_1,E) \to (X,0)$  be $R$ restricted to   $Y_1$  and  let $E=\rho {-1} (0) \subset Y_1$.  The total transform  of  $l_x\cup l_y$ by $\rho $, which is equal to  $E^+=\rho^{-1}(l_x\cup l_y)$, has a dual graph which   is a bamboo  as in Figure 5 with $k=m$  vertices. Let $l_x^1$ be the strict transform of $l_x$ by $\rho$. Then, $l_x^1$ only meets the  irreducible component of $E$ obtained by the last blowing-up of a line. The equation of $Y_1$ along $ l_x^1$ is $\{ z^r-xy^q=0\}$. 
  
 II)  If $r=1$, $Y_1$ is smooth and  Lemma 5.3.2 is proved i.e. $\rho_Y=\rho $. If $r>2$, the division $q=m'r+r'$ provides  $r'<r$. As $r$ is prime to $q$, $r'$ is prime to $r$ and $0<r'$. Moreover, we have  $r'<q$ because $r<q$. Let $R': Z' \to  Z$ be a sequence of $m'$ blowing-ups of the line $l_x ^1$ and of  its strict transforms. Let $Y_2$  be the strict transform of $Y_1$ by $R'$ and  let $\rho':(Y_2,E')\to (Y_1,E)$ be $R'$ restricted to $Y_2$.
    As $r<q$, $\rho' $ is bijective,  the dual graph of $\rho'^{-1}(E^+)$ is equal to the dual graph of $E^+$, which  is   a bamboo  as in Figure 5 with $k=m$  vertices. Moreover,  The equation of $Y_2$,  along the strict transform  of $ l_x^1$ by $\rho'$,   is $\{ z^r-xy^{r'}=0\}$. As $1\leq r'<r$ with $r$ and $r'$ prime to each other, Lemma 5.3.2 is proved by induction. 
    
    Let us justify the above statements I) and II)  by an  explicit computation of the blowing-up  of $l_x$.
  We consider  $Z_1=\{ ((x,y,z),(v:w)) \in \C ^3 \times \C P ^1, s.\ t.\ wy-vz=0 \}.$
  By definition,  the blowing-up  of $l_x$ in $\C^3 $, $R_1: Z_1 \to \C^3$,  is  the projection on $\C^3 $  restricted to $Z_1.$

As in  statement I), we consider   $X=\{(x,y,z)\in \C ^3  \ s.t.\  z^n-xy^q=0\}$ with $q<n$. We have to  describe  the strict transform $Y_{11}$ of $(X,0)$ by $R_1$,   the restriction $\rho_1 : (Y_{11},E) \to (X,0)$ of $R_1$ to $Y_{11}$,  $E_1=\rho_1^{-1}(0) $ and $E_1^+=\rho_1^{-1}(l_x\cup L_y) $. 

 i) In the chart $v=1$, we  have $(Z_1 \cap \{v=1\})=\{ ((x,y,wy),(1:w)) \in \C ^3 \times \C P ^1\}.$  The equation of $R_1^{-1}(0) \cap  \{v=1\}$  and of $E_1 \cap  \{v=1\} $   is $y=0$. 
    
    The equation of $(R_1^{-1}(X) \cap  \{v=1\})=(Y_{11}\cap  \{v=1\})$ is $\{ w^ny^{n-q}-x=0 \}$.   So,  all the points of $(\{v=1\}\cap Y_{11})$ are non singular and $( \{v=1\} \cap \{x\neq 0\} \cap Y_{11})$ doesn't meet $E_1$. 
    
    The strict transform of $l_x$ is not in $Y_{11} \cap \{v=1\}$. If $x=0$, we have: $$E_1\cap \{v=1\}=\{((0,0,0),(1:w)), w\in \C\} \subset Y_{11}.$$
    In $Y_{11}$, the strict  transform  $\tilde{l_y} =\{((0,y,0),(1:0)), y\in \C\} $ of $l_y$  meets $E_1$ at  $((0,0,0),(1:0))$.
    
   ii) In the chart $w=1$, we  have $(Z_1 \cap \{w=1\})=\{ ((x,vz, z),(v:1)) \in \C ^3 \times \C P ^1\}$.  The equation of $R_1^{-1}(0) \cap  \{w=1\}$  and of $E_1 \cap  \{w=1\} $   is $z=0$.

    The equation of  $(Y_{11}\cap  \{w=1\})$ is $\{ z^{n-q}-xv^q=0 \}$. So, the strict transform of $l_x$ is equal to
     $$ \tilde{l_x}=( \{w=1\}  \cap Y_{11} \cap R_1^{-1} (l_x))=\{ ((x,0, 0),(0:1)) \in \C ^3 \times \C P ^1\}.$$  The strict transform  $\tilde{l_x}$ meets $E_1$ at the point   $p_1=E_1 \cap \tilde{l_x}=((0,0,0),(0:1)$. 
     Then, $E_1= ((0,0,0) \times \C P^1)$  is included in $ Y_{11}$, moreover,  $ \tilde{l_x}$ and $ \tilde{l_y}$ meet $E_1$ at two distinct points. The total transform  $E_1^+=\rho_1^{-1}(l_x\cup L_y) $ consists in  one irreducible component $E_1$ and two germs of curves which meet $E_1$ in two distinct points. Moreover the equation of $Y_{11}$ along its singular locus $\tilde{l_x}$ is $\{ z^{n-q}-xy^q=0 \}$. By induction we obtain, as stated in I),  the germ  $(Y_1,0)$ defined by  $\{ z^r-xy^q=0\}$ with $1\leq r=n-mq<q$. 
     
 To justify  statement II), we  again consider  the blowing-up of $l_x$,  $\ R_1: Z_1 \to \C^3$. Let $Y_{12}$   be the strict transform of $Y_1$ by  $R_1$ and  let  $\rho_1': Y_{12} \to  Y_1$ be $R_1$ restricted to $Y_{12}$.  Then,   $Y_{12}$   has the equation $\{w^r-xy^{q-r}=0\}$ in the chart $v=1$. For all $x\in\C $, the intersection of $Y_{12}$ with $y=0$ is the only point $((x,0,0)),(1:0))$. In the chart $w=1$,  $Y_{12}$    has the equation $\{1-xv^q z^{q-r}=0\}$ and has empty intersection with $z=0$. This proves that $\rho'_1$ is bijective and by induction the map $\rho':(Y_2,E')\to (Y_1,E)$ describe above in II) is also bijective.

 {\it End of proof}.

{\bf Examples:}  

1) Let us consider $X=\{(x,y,z)\in \C ^3  \ s.t.\  z^n-xy^{n-1}=0\}$. The link of $(X,0)$ is the lens space L(n,1).  Let $R_1: Z_1 \to \C^3$ be the  blowing-up of the line $l_x$ in $\C^3$.  Let $Y$ be the strict transform of $X$ by $R_1$. The equation of $Y$ along the strict transform of $l_x$ is $\{ z-xy^{n-1}=0\}$. So,  $Y$ is non singular and we have obtained a resolution of $X$. Here the dual graph of the total  transform of $l_x\cup l_y$ is as in Figure 5 with only one vertex.

2) Let us consider $X=\{(x,y,z)\in \C ^3  \ s.t.\  z^n-xy^{n-2}=0\}$ with $n$ odd and $3<n$. The link of $(X,0)$ is the lens space L(n,2).  Let $R_1: Z_1 \to \C^3$ be the  blowing-up of  the line $l_x$ in $\C^3$. The equation of the strict transform  $Y_1$,  of $X$ by $R_1$,  along the strict transform of $l_x$ is $\{ z^2-xy^{n-2}=0\}$.  Let  $\rho : (Y_1,E) \to (X,0)$  be $R_1$ restricted to   $Y_1$. We write $n=2m+3$. As proved above, after $m $ blowing-ups  of lines, we obtain a surface $Y_2$ and a bijective morphism $\rho':(Y_2,E')\to (Y_1,E)$  such that the equation of $Y_2$ along the strict transform of $l_x$ is  $\{ z^2-xy=0\}$. The exceptional  divisor  $E$  of $\rho$ (resp. $E'$of $(\rho \circ \rho')$)  is an irreducible smooth rational curve. The blowing-up $\rho''$, of the intersection point   between $E'$ and the strict transform of $l_x$,  is a resolution of $Y_2$ and the exceptional divisor of $\rho''$ is a smooth rational curve. Then, $\rho \circ \rho' \circ \rho''$ is a resolution of $X=\{(x,y,z)\in \C ^3  \ s.t.\  z^n-xy^{n-2}=0\}$, the dual graph of its  exceptional divisor is a bamboo with two vertices.

\newpage

\section{ \bf An example of Hirzebruch-Jung's  resolution} 

We give the Hirzebruch-Jung  resolution  of the germ of surface in $\C ^3$ which satisfies  the following  equation: 
$$z^2= (x-y+y^3)(x-y+y^2)(y^{34}-(x-y)^{13}).$$

 where  $\pi : (X,0)\to (\C^2,0)$ is  the  projection on the $(x,y)$-plane.  It is a generic projection. In \cite{m-m} this example is  also detailed  for a non generic projection.

 The discriminant  locus of $\pi =(f,g)$ is the curve $\Delta $ which admits three components with Puiseux expansions given by :
  $$ x=y-y^2 $$
$$x= y-y^3$$
$$\ \ \   \ x=y+y^{34/13}$$
 
 Notice that the three components of $\Delta $ admit 1 as  first Puiseux exponent and respectively $2,3,34/13$ as second Puiseux exponent.
 
\begin{figure}[h]

$$
\unitlength=0.50mm
\begin{picture}(0.00,40.00)
\thicklines
\hspace*{10mm}
\put(-100,0){\line(1,0){160}}

\put(-80,0){\circle*{4}}
\put(-60,0){\circle*{4}}
\put(-40,0){\circle*{4}}
\put(-20,0){\circle*{4}}
\put(0,0){\circle*{4}}
\put(20,0){\circle*{4}}
\put(40,0){\circle*{4}}
\put(60,0){\circle*{4}}
\put(-60,0){\line(0,1){20}}
\put(-62,18){$* $}
\put(0,0){\line(0,1){20}}
\put(-2,18){$*$}
\put(58,18){$*$}
\put(60,0){\line(0,1){20}}

\put(-104,-2){$<$}

\put(-65,-10){$E_2$}
\put(-45,-10){$E_3$}
\put(-28,-10){$E_4$}
\put(-8,-10){$E_5$}
\put(15,-10){$E_6$}
\put(35,-10){$E_7$}
\put(55,-10){$E_8$}
\put(-85,-10){$E_1$}

\end{picture}
$$

\

\caption{The  dual graph of the minimal resolution of $\Delta $. An  irreducible component of the strict transform of $\Delta $  is  represented by  an  edge with a star. An edge  ended by  an  arrow represent the strict transform of $\{ x=0\}.$
} \label{figdisc3.1}
\end{figure}
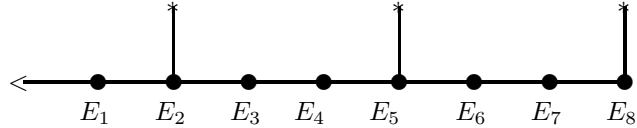

 The
 coordinate axes are transverse to the discriminant  locus of $\pi$. The dual graph $G(Z)$ is in Figure \ref{figdisc3.1}.

%

\

The dual  graph $G(\bar{Z})$ of $E_{\bar{Z}}$  admits a cycle created by the normalization. The irreducible component $E'_9$ of $E_Y$ is obtained by the resolution $\bar \rho$.
The irreducible components of the exceptional divisor associated to the vertices of  $G(\bar{Z})$  and $G(Y)$ have a  genus equal to zero. 

\begin{figure}[h]
 
$$
\unitlength=0.50mm
\begin{picture}(0.00,40.00) 
\thicklines
\hspace*{-20mm}

\put(-20,0){\line(1,0){20}}

\put(20,20){\line(1,0){20}}

\put(20,-20){\line(1,0){20}}

\put(0,0){\line(1,1){20}}
\put(0,0){\line(1,-1){20}}

\put(40,20){\line(1,-1){20}}
\put(40,-20){\line(1,1){20}}
\put(60,0){\line(1,0){20}}
\put(-20,0){\circle*{4}}
\put(0,0){\circle*{4}}
\put(20,20){\circle*{4}}
\put(20,-20){\circle*{4}}
\put(40,20){\circle*{4}}
\put(40,-20){\circle*{4}}
\put(60,0){\circle*{4}}
\put(80,0){\circle*{4}}
\put(100,0){\circle*{4}}
\put(120,0){\circle*{4}}

\put(120,-15){$E_8'$}
\put(100,-15){$E_7'$}
\put(-27,-15){$E_1'$}
\put(-7,-15){$E_2'$}
\put(17,30){$E_{3(1)}'$}
\put(17,-35){$E'_{3(2)}$}
\put(40,30){$E_{4(1)}'$}
\put(40,-35){$E'_{4(2)}$}
\put(57,-15){$E_5'$}
\put(77,-15){$E_6'$}
\put(140,-15){$E_9'$}

\put(80,0){\line(1,0){60}}
\put(140,0){\circle*{4}}


\end{picture}
$$
\vspace*{15mm}
\caption{The dual  graph $G(Y)$ of the Hirzebruch-Jung resolution associated to  $\pi$.}
\end{figure}
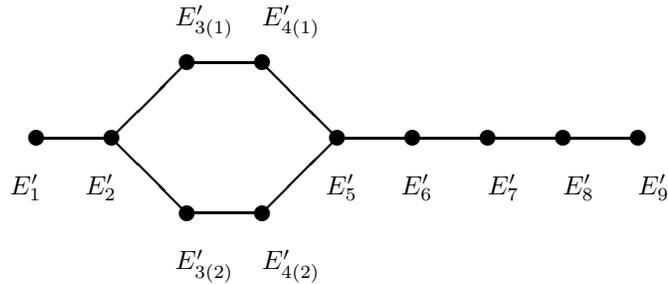

\begin{figure}[h]
 
$$
\unitlength=0.50mm
\begin{picture}(0.00,40.00)
\thicklines
\hspace*{-15mm}

\put(-20,0){\line(1,0){20}}
\put(20,20){\line(1,0){20}}
\put(20,-20){\line(1,0){20}}
\put(0,0){\line(1,1){20}}
\put(0,0){\line(1,-1){20}}
\put(40,20){\line(1,-1){20}}
\put(40,-20){\line(1,1){20}}
\put(60,0){\line(1,0){20}}
\put(-20,0){\circle*{4}}
\put(0,0){\circle*{4}}
\put(20,20){\circle*{4}}
\put(20,-20){\circle*{4}}
\put(40,20){\circle*{4}}
\put(40,-20){\circle*{4}}
\put(60,0){\circle*{4}}
\put(80,0){\circle*{4}}
\put(100,0){\circle*{4}}
\put(120,0){\circle*{4}}

\put(120,-15){$E_9'$}
\put(100,-15){$E_8'$}
\put(-27,-15){$E_1'$}
\put(-7,-15){$E_2'$}
\put(17,30){$E_{3(1)}'$}
\put(17,-35){$E'_{3(2)}$}
\put(40,30){$E_{4(1)}'$}
\put(40,-35){$E'_{4(2)}$}
\put(57,-15){$E_5'$}
\put(77,-15){$E_7'$}

\put(80,0){\line(1,0){40}}


\end{picture}
$$
\vspace*{15mm}
\caption{ The dual  graph $G(Y')$  of the minimal   resolution of $(X,0)$.}
\label{fig7}
\end{figure}
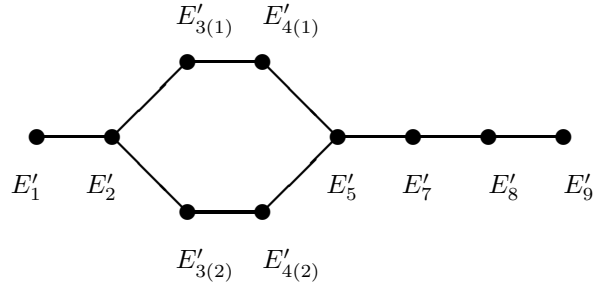

The  minimal good resolution $\rho$ is obtained by blowing down $E'_6$.  
Its dual graph is in Figure \ref{fig7}.

\vskip.1in

\noindent {\bf Adresses.}

\vskip.1in

\noindent Fran\c coise Michel / Laboratoire de Math\' ematiques Emile
Picard  /
Universit\' e Paul Sabatier / 118 route de Narbonne / F-31062
Toulouse / FRANCE

e-mail: fmichel@picard.ups-tlse.fr

\vskip.5in

\end{document}